\documentclass[11pt]{article}
\usepackage{graphicx}
\usepackage{amssymb,amsfonts,amsthm}

\usepackage{amsmath,amsthm,amssymb}
\usepackage{amscd}
\usepackage{amsfonts}
\usepackage{color}

\newcommand{\oss}[1]{}

\newcommand{\la}{\langle}
\newcommand{\ra}{\rangle}

\newcommand{\OO}{{\cal O}}

\newcommand{\iv}{^{-1}}

\newcommand{\F}{\mathbb{F}}

\newcommand{\HNN}{\hbox{HNN}}
\newcommand{\Q}{{\mathbb Q}}
\newcommand{\R}{{\mathbb R}}
\newcommand{\Z}{{\mathbb Z}}
\newcommand{\C}{{\mathbb C}}

\newcommand{\SL}{\hbox{SL}}

\newtheorem{theorem}{Theorem}[section]

\newtheorem{df}[theorem]{Definition}

\newtheorem{prob}[theorem]{Problem}

\title{Residual properties of 1-relator groups}
\author{Mark Sapir\thanks{This work was supported in part by the NSF grant DMS-0700811.}}
\date{}

\begin{document}

\maketitle

This text is based on my lectures given at the ``Groups St Andrews 2009" conference in Bath (August, 2009). I am going to survey the proof of the fact that almost all 1-related groups with at least 3 generators are residually finite (and satisfy some other nice properties) \cite{BS1, BS2, SS}, plus I will sketch proofs of some related facts. There are several unsolved problems at the end of the paper.

\section{Residually finite groups}

\begin{df} A group $G$ is called {\em residually finite} if for every $g\in G$, $g\ne 1$, there exists a homomorphism $\phi$ from $G$ onto a finite group $H$ such that $\phi(g)\ne 1$. If $H$ can be always chosen a $p$-group for some fixed prime $p$, then $G$ is called {\em residually (finite $p$-group).}
\end{df}

  {\bf Example.} Free groups $F_k$, cyclic extensions of free groups, finitely generated nilpotent groups, etc. are residually finite. The additive group $\Q$, infinite simple groups, free Burnside groups of sufficiently large exponents are not residually finite. The latter result follows from the combination of Novikov-Adyan's solution of the Bounded Burnside Problem \cite{NA} and Zelmanov's solution of the Restricted Burnside Problem \cite{Z1}.

\begin{theorem} Groups acting faithfully on rooted locally finite trees are residually finite.
\end{theorem}

\begin{center}
\unitlength .7mm 
\linethickness{0.4pt}
\ifx\plotpoint\undefined\newsavebox{\plotpoint}\fi 
\begin{picture}(116.5,54.75)(50,0)
\put(83.5,53.75){{\circle*{2}}}
\put(73.25,38.75){\circle*{2}} \put(91,39){\circle*{2}}
\multiput(83,54.25)(-.033737024,-.051903114){289}{\line(0,-1){.051903114}}
\multiput(82.75,54.25)(.03369565,-.06086957){230}{\line(0,-1){.06086957}}
\put(67,22.75){\circle*{2}} \put(84.75,23.25){\circle*{2}}
\put(78.25,23.25){\circle*{2}} \put(96,23.75){\circle*{2}}
\multiput(73,39)(-.03360215,-.08064516){186}{\line(0,-1){.08064516}}
\multiput(90.75,39.5)(-.03360215,-.08064516){186}{\line(0,-1){.08064516}}
\multiput(73.25,38.75)(.03368794,-.10460993){141}{\line(0,-1){.10460993}}
\multiput(91,39.25)(.03368794,-.10460993){141}{\line(0,-1){.10460993}}
\put(61.25,1.75){\circle*{1.5}} \put(72.75,2){\circle*{1.5}}
\put(90.75,2){\circle*{1.5}} \put(68.5,2){\circle*{1.58}}
\put(80,2.25){\circle*{1.58}} \put(98,2.25){\circle*{1.58}}
\multiput(67,22.5)(-.03360215,-.10887097){186}{\line(0,-1){.10887097}}
\multiput(78.5,22.75)(-.03360215,-.10887097){186}{\line(0,-1){.10887097}}
\multiput(96.5,22.75)(-.03360215,-.10887097){186}{\line(0,-1){.10887097}}
\multiput(67,23)(.0333333,-.4722222){45}{\line(0,-1){.4722222}}
\multiput(78.5,23.25)(.0333333,-.4722222){45}{\line(0,-1){.4722222}}
\multiput(96.5,23.25)(.0333333,-.4722222){45}{\line(0,-1){.4722222}}
\multiput(84.5,23.5)(-.0336538,-.3942308){52}{\line(0,-1){.3942308}}
\multiput(84.75,24)(.0333333,-.3541667){60}{\line(0,-1){.3541667}}
\put(82.5,2.75){\circle*{2}} \put(86.5,2.5){\circle*{2}}

\end{picture}

A rooted binary tree $T$.
\end{center}

\proof
Indeed, every automorphism $f$ of the tree $T$ must
fix the root and so it fixes the levels of the tree.
If $f\ne 1$ on level number $n$, we consider
the homomorphism from $\mathrm{Aut}(T)$ to the (finite) group of automorphisms of the finite tree consisting of the first $n$ levels of $T$. The homomorphisms
restricting automorphisms of $T$ to vertices of levels at most $n$.
The automorphism $f$ survives this homomorphism. Thus $\mathrm{Aut}(T)$ and all its subgroups are residually finite.
\endproof

Conversely (Kaluzhnin, see  \cite{GNS}) every finitely generated
residually finite group acts faithfully on a locally finite rooted
tree.

\subsection{Linear groups}

\begin{theorem}(A. Malcev, 1940, \cite{Mal}) Every finitely generated linear group is residually finite. Moreover, in characteristic 0, it is virtually residually (finite $p$-)group for all but finitely many primes $p$.
\end{theorem}

Note that a linear group itself may not be residually (finite $p$-)group for any $p$. Example: $\SL_3(\Z)$ by the Margulis' normal subgroup theorem \cite{Mar}.

\subsection{Problems}

Residual finiteness of a finitely presented group $G=\langle X \mid R\rangle$ is in fact a property of finite groups: we need to find out if there are ``enough" finite groups with $|X|$ generators satisfying the relations from $R$.

We shall consider two outstanding problems in group theory.

\begin{prob}\label{p1} (Gromov, \cite{Gr1}) Is every hyperbolic group
residually finite?
\end{prob}

This problem is hard because every non-elementary hyperbolic group has (by Gromov \cite{Gr1} and Olshanskii \cite{O1, O2, OlSq}) very many homomorphic images satisfying very strong finiteness properties (some of them are torsion and even of bounded torsion). This is in sharp contrast with CAT(0)-groups and even groups acting properly by isometries on finite dimensional CAT(0) cubical complexes which by a result of Wise \cite{Wise0} may be not residually finite and by a result of Burger and Mozes \cite{BM}  can be infinite and simple. By a result of Olshanskii \cite{OlBass}, Problem \ref{p1} has positive solution if and only if any two non-elementary hyperbolic groups have infinitely many common finite simple non-abelian factors.

\begin{prob}\label{p2} When is a one-relator group $\la
 X\mid R=1\ra$ residually finite?
\end{prob}

This problem is related to Problem \ref{p1} because it is easy to see that ``most" 1-relator groups satisfy small cancelation condition $C'(\lambda)$ (with arbitrary small $\lambda$) and so are hyperbolic. Note, though, that residual finiteness of small cancelation groups is as far from being proved (or disproved) as residual finiteness of all hyperbolic groups. Probably the strongest non-probabilistic result concerning Problem \ref{p2} was obtained by D. Wise \cite{Wise}: every 1-related group with a positive relator satisfying small cancelation condition $C'(\lambda)$, $\lambda \le 1/6$ is residually finite.

\bigskip

{\bf Example.} $BS(2,3)=\la a,
t\mid ta^2t\iv=a^3\ra$ is not residually finite ($a\mapsto a^2, t\mapsto t$ is a non-injective surjective endomophism showing that the group is not even Hopfian).

{\bf Example.} $BS(1,2) = \la a,t\mid
tat\iv =a^2\ra$ is metabelian, and representable by matrices over commutative ring (as all finitely generated metabelian groups), so it is residually finite by Malcev's result \cite{Mal}.

\subsection{The result}

The main goal of this article is to explain the proof  that the probabilistic version of Problem \ref{p2} has positive solution. Namely we are going to discuss the following

\begin{theorem}\label{thmain} (Borisov, Sapir, \v Spakulov\' a \cite{BS1}, \cite{BS2}, \cite{SS})
{Almost surely as $n\to \infty$, every 1-relator group with 3 or more generators and relator of length $n$, is}
\begin{itemize}
\item Residually finite,
\item Virtually residually (finite $p$-)group for all but finitely many primes $p$,
\item Coherent (that is all finitely generated subgroups are finitely presented).\end{itemize}
\end{theorem}

\subsection{Three probabilistic models}

The words ``almost surely" in Theorem \ref{thmain} need clarification. There are three natural probabilistic models to consider.

\begin{itemize}
\item {\bf Model 1.} Uniform distribution on words of length $\le n$.

\item {\bf Model 2.} Uniform distribution on cyclically reduced words of length $\le n$.

\item {\bf Model 3.} Uniform distribution on 1-relator groups given by cyclically reduced relators of length $\le n$ (up
to isomorphism)

\end{itemize}

These models turned out to be equivalent \cite{SS}. We prove that if the statement of the theorem is true for Model 2, then it is true for the other two models. Note that the equivalence Model $3\equiv$ Model $2$ uses a strong result of  Kapovich-Schupp-Shpilrain \cite{KSS}.

\section{Some properties of $1$-relator groups with 3 or more generators}

The following property has been mentioned above.

\bigskip

{\bf Fact 1. (Gromov \cite{Gr1})} {Almost every
1-relator group is hyperbolic.}

\bigskip
{\bf Fact 2. (Sacerdote and Schupp \cite{SSch})} {Every 1-relator group with 3 or more generators is $SQ$-universal (that is every countable group embeds into a quotient of that group.}

\bigskip

{\bf Fact 3 (B. Baumslag-Pride \cite{BP})}  Every group with the number of generators minus the number of relators  at least 2 is {\em large}, that is it has a subgroup of finite index that maps onto $F_2$.

\medskip
{Fact 3 and a result of P. Neumann \cite{N} imply Fact 2.}

Here is a close to complete proof of the result of Baumslag and Pride.

Let $G=\la x_1,...,x_g\mid u_1,...,u_r\ra$, $g-r\ge 2$.
First note that if $r-g\ge 1$ then we can assume that $t=x_1$ occurs in each relation with total exponent 0, i.e. $G$ maps onto $(Z,+)$, $t\to 1$. 

Rewrite the relators $u_k$ in terms of $s_{j,i}=t^ix_jt^{-i}$. Assume that for each $j\ge 2$, $m$ generators are involved: $0\le i\le m-1$.

Consider the subgroup $H$ that is a normal closure of $s=t^n$ and all $x_i, i\ge 2$, the kernel of the map $G\to \Z/n\Z$. It is of index $n$.

Consider the homomorphic image $\bar H$ of $H$ obtained by killing $s_{j,i}$, $2\le j\le g, 0\le i\le m-1$.

The standard Reidermeister-Schreier shows that $\bar H$ has $(g-1) (n-m) + 1$ generators $s, s_{j,i}$, $j=2,...,g, m\le i \le n$, and $nr$ relators not involving $s$.

So $\bar H$ is $\la s\ra * K$ where $K$ has $(g-1)(n-m)$ generators and $nr$ relators. For large enough $n$, then $\#$generators - $\#$relators of $K$ is $\ge 1$. So $K$ maps onto $\Z$, and $\bar H$ maps onto $F_2$. Q.E.D.

\subsection{Lackenby's result}

\begin{theorem} (Lackenby \cite{La}) For every large group $G$ and every $g\in G$ there exists $n$ such that $G/\la\la g^n\ra\ra$ is large.
\end{theorem}

The (very easy) proof given by Olshanskii and Osin \cite{OO} is essentially the same as the proof of the result of Baumslag and Pride above.

{{\bf Application.} There exists an infinite finitely generated group that is:}

\begin{itemize}
\item {residually finite
}
\item {torsion
}
\item {all sections are residually finite}

\item every finite section is solvable; every nilpotent finite section is Abelian.

\end{itemize}

\section{The Magnus procedure}
\bigskip

In order to deal with 1-relator groups, the main tool is the procedure invented by Magnus in the 30s. Here is an example.

{\bf Example (Magnus procedure).} {Consider the group}
$$\la a,b \mid aba\iv b\iv aba\iv b\iv a\iv b\iv
a=1\ra$$
For simplicity, we chose a relator with total exponent of $a$ equal 0 (as in the proof of Baumslag-Pride above, the general case reduces to this). We can write the relator as

$$(aba\iv)(b\iv)(aba\iv)(b\iv)(a\iv b\iv
a).$$

Replace $a^iba^{-i}$ by
$b_i$. The index $i$ is called {\em the Magnus $a$-index} of that letter. So we have
a new presentation of the same group.

$$
\la b_{-1}, b_0, b_1, a \mid b_1b_0\iv b_1b_0\iv b_{-1}=1, ab_{-1}a\iv = b_0, ab_0a\iv=b_1\ra.$$
Note that $b_{-1}$
appears only once in $b_1b_0\iv b_1b_0\iv b_{-1}\iv$.
So we can replace
$b_{-1}$ by $b_1b_0\iv b_1b_0\iv$, remove this generator, and get a
new presentation of the same group.

$$\la b_0, b_1, a \mid a\iv b_0 a=b_1b_0\iv b_1b_0\iv, a\iv b_1a=b_0\ra.$$

This is clearly an ascending HNN
extension of the free group {$\la b_0,b_1\ra$}. Thus the initial group is an ascending HNN extension of a free group. We shall see below that this happens quite often.

\subsection{Ascending HNN extensions}

{\bf Definition.} Let $G$ be a group, $\phi\colon G\to G$ be an
injective endomorphism.   The group
{$$\HNN_\phi(G)=\la G, t\mid tat\iv=\phi(a), a\in
G\ra$$}is called an {\em ascending HNN extension} of $G$ or the {\em
mapping torus} of $\phi$.

  {\bf Example.} Here is the main motivational example for us: $H_{T}=\la x, y, t\mid txt\iv = xy,
tyt\iv=yx\ra$ corresponding to the Thue endomorphism of the free group. We shall return to this example several times later.

Ascending HNN extensions have nice geometric interpretations. They are fundamental groups of the mapping tori of the endomorphisms. On the picture below, the handle is attached on one side to the curves $x$ and $y$, and on the other side to $xy$ and $yx$. The fundamental group of this object is $H_T$.

\unitlength .5mm 
\linethickness{0.4pt}
\ifx\plotpoint\undefined\newsavebox{\plotpoint}\fi 
\begin{picture}(170.5,97.88)(0,0)
\put(58.99,62.25){\line(0,1){.728}}
\put(58.98,62.98){\line(0,1){.726}}
\put(58.92,63.7){\line(0,1){.723}}
\multiput(58.83,64.43)(-.0319,.1793){4}{\line(0,1){.1793}}
\multiput(58.7,65.14)(-.03268,.14199){5}{\line(0,1){.14199}}
\multiput(58.54,65.85)(-.03315,.11681){6}{\line(0,1){.11681}}
\multiput(58.34,66.56)(-.03341,.09857){7}{\line(0,1){.09857}}
\multiput(58.11,67.25)(-.03353,.08467){8}{\line(0,1){.08467}}
\multiput(57.84,67.92)(-.03355,.07367){9}{\line(0,1){.07367}}
\multiput(57.54,68.59)(-.03348,.0647){10}{\line(0,1){.0647}}
\multiput(57.2,69.23)(-.03336,.05722){11}{\line(0,1){.05722}}
\multiput(56.83,69.86)(-.03317,.05084){12}{\line(0,1){.05084}}
\multiput(56.44,70.47)(-.032943,.045335){13}{\line(0,1){.045335}}
\multiput(56.01,71.06)(-.032666,.040507){14}{\line(0,1){.040507}}
\multiput(55.55,71.63)(-.03235,.036227){15}{\line(0,1){.036227}}
\multiput(55.06,72.17)(-.031996,.032396){16}{\line(0,1){.032396}}
\multiput(54.55,72.69)(-.035823,.032797){15}{\line(-1,0){.035823}}
\multiput(54.02,73.18)(-.040098,.033166){14}{\line(-1,0){.040098}}
\multiput(53.45,73.65)(-.044923,.033502){13}{\line(-1,0){.044923}}
\multiput(52.87,74.08)(-.04655,.031203){13}{\line(-1,0){.04655}}
\multiput(52.26,74.49)(-.05206,.03123){12}{\line(-1,0){.05206}}
\multiput(51.64,74.86)(-.05844,.03117){11}{\line(-1,0){.05844}}
\multiput(51,75.21)(-.06592,.03101){10}{\line(-1,0){.06592}}
\multiput(50.34,75.52)(-.07489,.03073){9}{\line(-1,0){.07489}}
\multiput(49.66,75.79)(-.08588,.0303){8}{\line(-1,0){.08588}}
\multiput(48.98,76.03)(-.09976,.02965){7}{\line(-1,0){.09976}}
\multiput(48.28,76.24)(-.11798,.0287){6}{\line(-1,0){.11798}}
\multiput(47.57,76.41)(-.14313,.02728){5}{\line(-1,0){.14313}}
\put(46.86,76.55){\line(-1,0){.722}}
\put(46.13,76.65){\line(-1,0){.726}}
\put(45.41,76.72){\line(-1,0){.728}}
\put(44.68,76.74){\line(-1,0){.728}}
\put(43.95,76.73){\line(-1,0){.727}}
\put(43.22,76.69){\line(-1,0){.724}}
\put(42.5,76.6){\line(-1,0){.719}}
\multiput(41.78,76.49)(-.14239,-.03092){5}{\line(-1,0){.14239}}
\multiput(41.07,76.33)(-.11721,-.03169){6}{\line(-1,0){.11721}}
\multiput(40.37,76.14)(-.09897,-.03218){7}{\line(-1,0){.09897}}
\multiput(39.67,75.92)(-.08508,-.03247){8}{\line(-1,0){.08508}}
\multiput(38.99,75.66)(-.07408,-.03263){9}{\line(-1,0){.07408}}
\multiput(38.33,75.36)(-.06511,-.03268){10}{\line(-1,0){.06511}}
\multiput(37.68,75.04)(-.05762,-.03265){11}{\line(-1,0){.05762}}
\multiput(37.04,74.68)(-.05125,-.03254){12}{\line(-1,0){.05125}}
\multiput(36.43,74.29)(-.045741,-.032378){13}{\line(-1,0){.045741}}
\multiput(35.83,73.87)(-.040909,-.032161){14}{\line(-1,0){.040909}}
\multiput(35.26,73.42)(-.036625,-.031898){15}{\line(-1,0){.036625}}
\multiput(34.71,72.94)(-.034976,-.033698){15}{\line(-1,0){.034976}}
\multiput(34.18,72.43)(-.033238,-.035413){15}{\line(0,-1){.035413}}
\multiput(33.69,71.9)(-.033661,-.039684){14}{\line(0,-1){.039684}}
\multiput(33.22,71.34)(-.031624,-.041325){14}{\line(0,-1){.041325}}
\multiput(32.77,70.77)(-.031778,-.046159){13}{\line(0,-1){.046159}}
\multiput(32.36,70.17)(-.03187,-.05167){12}{\line(0,-1){.05167}}
\multiput(31.98,69.55)(-.03189,-.05805){11}{\line(0,-1){.05805}}
\multiput(31.63,68.91)(-.03183,-.06553){10}{\line(0,-1){.06553}}
\multiput(31.31,68.25)(-.03166,-.0745){9}{\line(0,-1){.0745}}
\multiput(31.02,67.58)(-.03136,-.0855){8}{\line(0,-1){.0855}}
\multiput(30.77,66.9)(-.03089,-.09939){7}{\line(0,-1){.09939}}
\multiput(30.56,66.2)(-.03016,-.11761){6}{\line(0,-1){.11761}}
\multiput(30.37,65.5)(-.02906,-.14278){5}{\line(0,-1){.14278}}
\put(30.23,64.78){\line(0,-1){.72}}
\put(30.12,64.06){\line(0,-1){.725}}
\put(30.05,63.34){\line(0,-1){2.909}}
\put(30.12,60.43){\line(0,-1){.72}}
\multiput(30.23,59.71)(.02915,-.14276){5}{\line(0,-1){.14276}}
\multiput(30.38,58.99)(.03024,-.1176){6}{\line(0,-1){.1176}}
\multiput(30.56,58.29)(.03095,-.09937){7}{\line(0,-1){.09937}}
\multiput(30.77,57.59)(.03142,-.08547){8}{\line(0,-1){.08547}}
\multiput(31.03,56.91)(.03171,-.07448){9}{\line(0,-1){.07448}}
\multiput(31.31,56.24)(.03187,-.06551){10}{\line(0,-1){.06551}}
\multiput(31.63,55.58)(.03193,-.05803){11}{\line(0,-1){.05803}}
\multiput(31.98,54.95)(.0319,-.05165){12}{\line(0,-1){.05165}}
\multiput(32.36,54.33)(.031808,-.046139){13}{\line(0,-1){.046139}}
\multiput(32.78,53.73)(.031651,-.041305){14}{\line(0,-1){.041305}}
\multiput(33.22,53.15)(.033687,-.039662){14}{\line(0,-1){.039662}}
\multiput(33.69,52.59)(.033261,-.035391){15}{\line(0,-1){.035391}}
\multiput(34.19,52.06)(.034998,-.033675){15}{\line(1,0){.034998}}
\multiput(34.72,51.56)(.036646,-.031874){15}{\line(1,0){.036646}}
\multiput(35.27,51.08)(.04093,-.032135){14}{\line(1,0){.04093}}
\multiput(35.84,50.63)(.045762,-.032348){13}{\line(1,0){.045762}}
\multiput(36.43,50.21)(.05127,-.03251){12}{\line(1,0){.05127}}
\multiput(37.05,49.82)(.05765,-.03261){11}{\line(1,0){.05765}}
\multiput(37.68,49.46)(.06513,-.03264){10}{\line(1,0){.06513}}
\multiput(38.33,49.13)(.0741,-.03258){9}{\line(1,0){.0741}}
\multiput(39,48.84)(.0851,-.03242){8}{\line(1,0){.0851}}
\multiput(39.68,48.58)(.099,-.03212){7}{\line(1,0){.099}}
\multiput(40.38,48.36)(.11723,-.03162){6}{\line(1,0){.11723}}
\multiput(41.08,48.17)(.14241,-.03082){5}{\line(1,0){.14241}}
\put(41.79,48.01){\line(1,0){.719}}
\put(42.51,47.89){\line(1,0){.724}}
\put(43.23,47.81){\line(1,0){.727}}
\put(43.96,47.77){\line(1,0){.728}}
\put(44.69,47.76){\line(1,0){.728}}
\put(45.42,47.79){\line(1,0){.726}}
\put(46.14,47.85){\line(1,0){.722}}
\multiput(46.86,47.95)(.14311,.02737){5}{\line(1,0){.14311}}
\multiput(47.58,48.09)(.11796,.02878){6}{\line(1,0){.11796}}
\multiput(48.29,48.26)(.09974,.02971){7}{\line(1,0){.09974}}
\multiput(48.99,48.47)(.08586,.03035){8}{\line(1,0){.08586}}
\multiput(49.67,48.71)(.07487,.03078){9}{\line(1,0){.07487}}
\multiput(50.35,48.99)(.0659,.03105){10}{\line(1,0){.0659}}
\multiput(51.01,49.3)(.05842,.03121){11}{\line(1,0){.05842}}
\multiput(51.65,49.64)(.05204,.03126){12}{\line(1,0){.05204}}
\multiput(52.27,50.02)(.04653,.031233){13}{\line(1,0){.04653}}
\multiput(52.88,50.42)(.044901,.033532){13}{\line(1,0){.044901}}
\multiput(53.46,50.86)(.040077,.033192){14}{\line(1,0){.040077}}
\multiput(54.02,51.32)(.035801,.03282){15}{\line(1,0){.035801}}
\multiput(54.56,51.82)(.031975,.032416){16}{\line(0,1){.032416}}
\multiput(55.07,52.33)(.032326,.036248){15}{\line(0,1){.036248}}
\multiput(55.56,52.88)(.03264,.040528){14}{\line(0,1){.040528}}
\multiput(56.01,53.45)(.032913,.045357){13}{\line(0,1){.045357}}
\multiput(56.44,54.04)(.03314,.05087){12}{\line(0,1){.05087}}
\multiput(56.84,54.65)(.03332,.05724){11}{\line(0,1){.05724}}
\multiput(57.2,55.28)(.03344,.06472){10}{\line(0,1){.06472}}
\multiput(57.54,55.92)(.0335,.07369){9}{\line(0,1){.07369}}
\multiput(57.84,56.59)(.03347,.08469){8}{\line(0,1){.08469}}
\multiput(58.11,57.26)(.03334,.09859){7}{\line(0,1){.09859}}
\multiput(58.34,57.95)(.03307,.11683){6}{\line(0,1){.11683}}
\multiput(58.54,58.65)(.03259,.14201){5}{\line(0,1){.14201}}
\multiput(58.7,59.36)(.0318,.1793){4}{\line(0,1){.1793}}
\put(58.83,60.08){\line(0,1){.723}}
\put(58.92,60.8){\line(0,1){1.445}}
\put(57.74,33.5){\line(0,1){.728}}
\put(57.73,34.23){\line(0,1){.726}}
\put(57.67,34.95){\line(0,1){.723}}
\multiput(57.58,35.68)(-.0319,.1793){4}{\line(0,1){.1793}}
\multiput(57.45,36.39)(-.03268,.14199){5}{\line(0,1){.14199}}
\multiput(57.29,37.1)(-.03315,.11681){6}{\line(0,1){.11681}}
\multiput(57.09,37.81)(-.03341,.09857){7}{\line(0,1){.09857}}
\multiput(56.86,38.5)(-.03353,.08467){8}{\line(0,1){.08467}}
\multiput(56.59,39.17)(-.03355,.07367){9}{\line(0,1){.07367}}
\multiput(56.29,39.84)(-.03348,.0647){10}{\line(0,1){.0647}}
\multiput(55.95,40.48)(-.03336,.05722){11}{\line(0,1){.05722}}
\multiput(55.58,41.11)(-.03317,.05084){12}{\line(0,1){.05084}}
\multiput(55.19,41.72)(-.032943,.045335){13}{\line(0,1){.045335}}
\multiput(54.76,42.31)(-.032666,.040507){14}{\line(0,1){.040507}}
\multiput(54.3,42.88)(-.03235,.036227){15}{\line(0,1){.036227}}
\multiput(53.81,43.42)(-.031996,.032396){16}{\line(0,1){.032396}}
\multiput(53.3,43.94)(-.035823,.032797){15}{\line(-1,0){.035823}}
\multiput(52.77,44.43)(-.040098,.033166){14}{\line(-1,0){.040098}}
\multiput(52.2,44.9)(-.044923,.033502){13}{\line(-1,0){.044923}}
\multiput(51.62,45.33)(-.04655,.031203){13}{\line(-1,0){.04655}}
\multiput(51.01,45.74)(-.05206,.03123){12}{\line(-1,0){.05206}}
\multiput(50.39,46.11)(-.05844,.03117){11}{\line(-1,0){.05844}}
\multiput(49.75,46.46)(-.06592,.03101){10}{\line(-1,0){.06592}}
\multiput(49.09,46.77)(-.07489,.03073){9}{\line(-1,0){.07489}}
\multiput(48.41,47.04)(-.08588,.0303){8}{\line(-1,0){.08588}}
\multiput(47.73,47.28)(-.09976,.02965){7}{\line(-1,0){.09976}}
\multiput(47.03,47.49)(-.11798,.0287){6}{\line(-1,0){.11798}}
\multiput(46.32,47.66)(-.14313,.02728){5}{\line(-1,0){.14313}}
\put(45.61,47.8){\line(-1,0){.722}}
\put(44.88,47.9){\line(-1,0){.726}}
\put(44.16,47.97){\line(-1,0){.728}}
\put(43.43,47.99){\line(-1,0){.728}}
\put(42.7,47.98){\line(-1,0){.727}}
\put(41.97,47.94){\line(-1,0){.724}}
\put(41.25,47.85){\line(-1,0){.719}}
\multiput(40.53,47.74)(-.14239,-.03092){5}{\line(-1,0){.14239}}
\multiput(39.82,47.58)(-.11721,-.03169){6}{\line(-1,0){.11721}}
\multiput(39.12,47.39)(-.09897,-.03218){7}{\line(-1,0){.09897}}
\multiput(38.42,47.17)(-.08508,-.03247){8}{\line(-1,0){.08508}}
\multiput(37.74,46.91)(-.07408,-.03263){9}{\line(-1,0){.07408}}
\multiput(37.08,46.61)(-.06511,-.03268){10}{\line(-1,0){.06511}}
\multiput(36.43,46.29)(-.05762,-.03265){11}{\line(-1,0){.05762}}
\multiput(35.79,45.93)(-.05125,-.03254){12}{\line(-1,0){.05125}}
\multiput(35.18,45.54)(-.045741,-.032378){13}{\line(-1,0){.045741}}
\multiput(34.58,45.12)(-.040909,-.032161){14}{\line(-1,0){.040909}}
\multiput(34.01,44.67)(-.036625,-.031898){15}{\line(-1,0){.036625}}
\multiput(33.46,44.19)(-.034976,-.033698){15}{\line(-1,0){.034976}}
\multiput(32.93,43.68)(-.033238,-.035413){15}{\line(0,-1){.035413}}
\multiput(32.44,43.15)(-.033661,-.039684){14}{\line(0,-1){.039684}}
\multiput(31.97,42.59)(-.031624,-.041325){14}{\line(0,-1){.041325}}
\multiput(31.52,42.02)(-.031778,-.046159){13}{\line(0,-1){.046159}}
\multiput(31.11,41.42)(-.03187,-.05167){12}{\line(0,-1){.05167}}
\multiput(30.73,40.8)(-.03189,-.05805){11}{\line(0,-1){.05805}}
\multiput(30.38,40.16)(-.03183,-.06553){10}{\line(0,-1){.06553}}
\multiput(30.06,39.5)(-.03166,-.0745){9}{\line(0,-1){.0745}}
\multiput(29.77,38.83)(-.03136,-.0855){8}{\line(0,-1){.0855}}
\multiput(29.52,38.15)(-.03089,-.09939){7}{\line(0,-1){.09939}}
\multiput(29.31,37.45)(-.03016,-.11761){6}{\line(0,-1){.11761}}
\multiput(29.12,36.75)(-.02906,-.14278){5}{\line(0,-1){.14278}}
\put(28.98,36.03){\line(0,-1){.72}}
\put(28.87,35.31){\line(0,-1){.725}}
\put(28.8,34.59){\line(0,-1){2.909}}
\put(28.87,31.68){\line(0,-1){.72}}
\multiput(28.98,30.96)(.02915,-.14276){5}{\line(0,-1){.14276}}
\multiput(29.13,30.24)(.03024,-.1176){6}{\line(0,-1){.1176}}
\multiput(29.31,29.54)(.03095,-.09937){7}{\line(0,-1){.09937}}
\multiput(29.52,28.84)(.03142,-.08547){8}{\line(0,-1){.08547}}
\multiput(29.78,28.16)(.03171,-.07448){9}{\line(0,-1){.07448}}
\multiput(30.06,27.49)(.03187,-.06551){10}{\line(0,-1){.06551}}
\multiput(30.38,26.83)(.03193,-.05803){11}{\line(0,-1){.05803}}
\multiput(30.73,26.2)(.0319,-.05165){12}{\line(0,-1){.05165}}
\multiput(31.11,25.58)(.031808,-.046139){13}{\line(0,-1){.046139}}
\multiput(31.53,24.98)(.031651,-.041305){14}{\line(0,-1){.041305}}
\multiput(31.97,24.4)(.033687,-.039662){14}{\line(0,-1){.039662}}
\multiput(32.44,23.84)(.033261,-.035391){15}{\line(0,-1){.035391}}
\multiput(32.94,23.31)(.034998,-.033675){15}{\line(1,0){.034998}}
\multiput(33.47,22.81)(.036646,-.031874){15}{\line(1,0){.036646}}
\multiput(34.02,22.33)(.04093,-.032135){14}{\line(1,0){.04093}}
\multiput(34.59,21.88)(.045762,-.032348){13}{\line(1,0){.045762}}
\multiput(35.18,21.46)(.05127,-.03251){12}{\line(1,0){.05127}}
\multiput(35.8,21.07)(.05765,-.03261){11}{\line(1,0){.05765}}
\multiput(36.43,20.71)(.06513,-.03264){10}{\line(1,0){.06513}}
\multiput(37.08,20.38)(.0741,-.03258){9}{\line(1,0){.0741}}
\multiput(37.75,20.09)(.0851,-.03242){8}{\line(1,0){.0851}}
\multiput(38.43,19.83)(.099,-.03212){7}{\line(1,0){.099}}
\multiput(39.13,19.61)(.11723,-.03162){6}{\line(1,0){.11723}}
\multiput(39.83,19.42)(.14241,-.03082){5}{\line(1,0){.14241}}
\put(40.54,19.26){\line(1,0){.719}}
\put(41.26,19.14){\line(1,0){.724}}
\put(41.98,19.06){\line(1,0){.727}}
\put(42.71,19.02){\line(1,0){.728}}
\put(43.44,19.01){\line(1,0){.728}}
\put(44.17,19.04){\line(1,0){.726}}
\put(44.89,19.1){\line(1,0){.722}}
\multiput(45.61,19.2)(.14311,.02737){5}{\line(1,0){.14311}}
\multiput(46.33,19.34)(.11796,.02878){6}{\line(1,0){.11796}}
\multiput(47.04,19.51)(.09974,.02971){7}{\line(1,0){.09974}}
\multiput(47.74,19.72)(.08586,.03035){8}{\line(1,0){.08586}}
\multiput(48.42,19.96)(.07487,.03078){9}{\line(1,0){.07487}}
\multiput(49.1,20.24)(.0659,.03105){10}{\line(1,0){.0659}}
\multiput(49.76,20.55)(.05842,.03121){11}{\line(1,0){.05842}}
\multiput(50.4,20.89)(.05204,.03126){12}{\line(1,0){.05204}}
\multiput(51.02,21.27)(.04653,.031233){13}{\line(1,0){.04653}}
\multiput(51.63,21.67)(.044901,.033532){13}{\line(1,0){.044901}}
\multiput(52.21,22.11)(.040077,.033192){14}{\line(1,0){.040077}}
\multiput(52.77,22.57)(.035801,.03282){15}{\line(1,0){.035801}}
\multiput(53.31,23.07)(.031975,.032416){16}{\line(0,1){.032416}}
\multiput(53.82,23.58)(.032326,.036248){15}{\line(0,1){.036248}}
\multiput(54.31,24.13)(.03264,.040528){14}{\line(0,1){.040528}}
\multiput(54.76,24.7)(.032913,.045357){13}{\line(0,1){.045357}}
\multiput(55.19,25.29)(.03314,.05087){12}{\line(0,1){.05087}}
\multiput(55.59,25.9)(.03332,.05724){11}{\line(0,1){.05724}}
\multiput(55.95,26.53)(.03344,.06472){10}{\line(0,1){.06472}}
\multiput(56.29,27.17)(.0335,.07369){9}{\line(0,1){.07369}}
\multiput(56.59,27.84)(.03347,.08469){8}{\line(0,1){.08469}}
\multiput(56.86,28.51)(.03334,.09859){7}{\line(0,1){.09859}}
\multiput(57.09,29.2)(.03307,.11683){6}{\line(0,1){.11683}}
\multiput(57.29,29.9)(.03259,.14201){5}{\line(0,1){.14201}}
\multiput(57.45,30.61)(.0318,.1793){4}{\line(0,1){.1793}}
\put(57.58,31.33){\line(0,1){.723}}
\put(57.67,32.05){\line(0,1){1.445}}
\put(43.5,47.75){\circle*{2.55}}

\qbezier(76.75,61.75)(121.88,97.88)(162.5,52.5)
\qbezier(77.5,42)(122.63,78.13)(163.25,32.75)
\qbezier(77.25,42)(122.38,78.13)(163,32.75)
\qbezier(78,22.25)(123.13,58.38)(163.75,13)
\qbezier(77.25,42.25)(67,53.63)(76.75,61.5)
\qbezier(78,22.5)(67.75,33.88)(77.5,41.75)
\qbezier(163.25,33.25)(153,44.63)(162.75,52.5)
\qbezier(164,13.5)(153.75,24.88)(163.5,32.75)
\multiput(76.75,61.5)(-.03125,-.03125){8}{\line(0,-1){.03125}}
\multiput(77.5,41.75)(-.03125,-.03125){8}{\line(0,-1){.03125}}
\multiput(76.5,61.25)(.033333,-.05){30}{\line(0,-1){.05}}
\multiput(77.25,41.5)(.033333,-.05){30}{\line(0,-1){.05}}
\multiput(78.25,59)(.032609,-.086957){23}{\line(0,-1){.086957}}
\multiput(79,39.25)(.032609,-.086957){23}{\line(0,-1){.086957}}
\multiput(79.25,55.5)(.03125,-.28125){8}{\line(0,-1){.28125}}
\multiput(80,35.75)(.03125,-.28125){8}{\line(0,-1){.28125}}
\multiput(79.5,51.5)(-.03125,-.34375){8}{\line(0,-1){.34375}}
\multiput(80.25,31.75)(-.03125,-.34375){8}{\line(0,-1){.34375}}
\multiput(78.5,47)(-.0328947,-.1118421){38}{\line(0,-1){.1118421}}
\multiput(79.25,27.25)(-.0328947,-.1118421){38}{\line(0,-1){.1118421}}
\qbezier(162.5,52.25)(169.75,50.25)(163,33.25)
\qbezier(163.25,32.5)(170.5,30.5)(163.75,13.5)

\put(47.25,67.5){\makebox(0,0)[cc]{$x$}}
\put(44.75,24){\makebox(0,0)[cc]{$y$}}

 \put(81.25,59.75){\makebox(0,0)[cc]{$x$}}
\put(82.25,35){\makebox(0,0)[cc]{$y$}}
\put(169,45.75){\makebox(0,0)[cc]{$xy$}}
\put(169.25,24){\makebox(0,0)[cc]{$yx$}}
\end{picture}

\subsection{Facts about ascending HNN extensions}

\begin{itemize}

\item {Every element in an ascending $\HNN$ extension of $G$ can be
represented in the form $t^{-k}gt^\ell$ for some $k,\ell\in
\mathbb{Z}$ and $g\in G$. The number $\ell-k$ is an invariant, the
representation is unique for a given $k$.}

\item {(Feighn-Handel \cite{FH}) If $G$ is free
then $\HNN_\phi(G)$ is {\em coherent} i.e. every f.g. subgroup is
f.p.}

\item {(Geoghegan-Mihalik-Sapir-Wise\cite{GMSW}) If
$G$ is free then $\HNN_\phi(G)$ is {\em Hopfian} i.e. every
surjective endomorphism is injective.}

\item {(Sapir-Wise\cite{SW} ) An ascending \HNN extension of a residually finite group can be non-residually finite (example - Grigorcuk's group and its Lysenok extension).}
\end{itemize}

\subsection{Walks in $\Z^2$}

Consider a 1-relator group $\la a,b \mid R\ra$ with 2 generators. We can consider the 2-dimensional grid $\Z^2$, label horizontal edges by $a$, vertical edges by $b$. Then the word $R$ corresponds to a walk in $\Z^2$ starting at $(0,0)$. Suppose that the total sum of exponents of $a$ in $R$ is 0.  Then the projections of the vertical steps of the walk onto the horizontal axis give the Magnus $a$-indices of the corresponding $b$'s. Hence the group is an ascending HNN extension of a free group if one of the two support vertical lines intersects the walk in exactly one step.

A similar fact holds when the total sum of exponents of $a$ is not 0.

\subsection{Results of Ken Brown}

\begin{theorem} (Ken Brown, \cite{Brown}) Let $G=\langle x_1,...,x_k\mid R=1\rangle$ be a 1-relator group. Let $w$ be the corresponding walk in $\Z^k$, connecting point $O$ with point $M$.
\begin{itemize}
\item If $k=2$. Then $G$ is an ascending HNN extension of a free group if and only if one of the two support lines  of $w$ that is parallel to $\vec{OM}$ (i.e. one of the two lines that are parallen to $\vec{OM}$, intersect $w$ and such that the whole $w$ is on one side of that line) intersects $w$ in a single vertex or a single edge.

    \item If $k>2$ then $G$ is never an ascending HNN extension of a free group.\end{itemize}
\end{theorem}

\subsection{The Dunfield-Thurston result}

Let $k=2$. What is the probability that Brown's condition holds?
Consider a random word $R$ in $a$ and $b$ and the corresponding random walk $w$ in $\Z^2$. The projection of $w$ onto the line $R$ perpendicular to $\vec{OM}$ is a random {\em bridge} (i.e walk that starts and ends at the same point). Support lines map to the two extreme points of the bridge.

A bridge is called {\em good} if it visits one of its extreme point only once, otherwise it is {\em bad}. Good bridges correspond to words which satisfy Brown's condition, bad bridges correspond to words that do not satisfy that condition. The number of good (bad) walks of length $n$ is denoted by $\# good(n)$ (resp. $\# bad(n)$). The probability for the walk to be good (bad) is denoted by $p_{good}$ (resp. $p_{bad}=1-p_{good}$).

Note that we can turn any good bridge into a bad one by inserting one or two subwords of length 8 in the corresponding word. Hence

$$\# good(n)\le \# bad(n+16).$$

The number of words of length $n+16$ is at most $4^{16}$ times the number of words of length $n$. Hence

\medskip

$$p_{good}\le 4^{16} p_{bad}$$

\medskip

Hence $p_{good}<1$. Similarly $p_{bad}<1$. In fact the Monte Carlo experiments (pick a large $n$, and a large number of group words $R$ in two generators; for each word $R$ check if Magnus rewriting of $R$ produces a word with the maximal or minimal index occurring only once) show that $p_{good}$ is about $.96$.

\subsection{The Congruence Extension Property}

 \begin{theorem} \label{O1} (Olshanskii \cite{OlSq} Let $K$ be a collection of (cyclic) words in $\{a,b\}$ that satisfy $C'(1/12)$. Then the subgroup $N$ of $F_2$ generated by $K$ satisfies the {\em congruence extension property} that is for every normal subgroup $L\lhd N$, $\la\la L\ra\ra_F\cap N=L$. Hence $H=N/L$ embeds into $G=F_2/\la\la L\ra\ra$.

 \end{theorem}

{\bf Remark.} In fact \cite{OlSq} contains a much stronger result for arbitrary hyperbolic groups.

\proof Consider $L$ as the set of relations of $G$. We need to show that the kernel of the natural map $N\to G$ is $L$, that is if $w(K)=1$ modulo $L$ (here $w(K)$ is obtained from $w$ by plugging elements of $K$ for its letters), then $w\in L$. Consider a van Kampen diagram $\Delta$ for the equality $w(K)=1$ with minimal possible number of cells. The boundary of every cell is a product of words from $K$, called {\em blocks}. If two cells touch by a subpath of their boundaries that includes a ``large" portion (say, $\frac{1}{12}$) of a block, then by the small cancelation condition, these boundaries have a common block. Hence we can remove the common part of the boundary consisting of several blocks, and form a new cell with boundary label still in $L$. That contradicts the assumption that $\Delta$ is minimal. This implies that cells in $\Delta$ do not have common large parts of their boundaries. Hence $\Delta$ is a small cancelation map. By the standard Greendlinger lemma, one of the cells $\pi$ in $\Delta$ has a large part (say, $> \frac{1}{1/12}$) of a block in common with the boundary of $\Delta$. Then we can cut $\pi$ out of $\Delta$ and obtain a smaller diagram with the same property. The proofs concludes by induction. The proof is illustrated by the following picture.  The left part of it illustrates ``cancelation" of cells with large common parts of the boundaries. The right part illustrates the induction step.

\unitlength .5mm 
\linethickness{0.4pt}
\ifx\plotpoint\undefined\newsavebox{\plotpoint}\fi 
\begin{picture}(194.875,95.5)(0,0)
\put(31.146,21.5){\line(0,1){.7245}}
\put(31.128,22.225){\line(0,1){.7227}}
\put(31.073,22.947){\line(0,1){.719}}
\multiput(30.982,23.666)(-.03176,.17838){4}{\line(0,1){.17838}}
\multiput(30.855,24.38)(-.032562,.141247){5}{\line(0,1){.141247}}
\multiput(30.692,25.086)(-.033024,.116191){6}{\line(0,1){.116191}}
\multiput(30.494,25.783)(-.033283,.098042){7}{\line(0,1){.098042}}
\multiput(30.261,26.469)(-.033403,.084212){8}{\line(0,1){.084212}}
\multiput(29.994,27.143)(-.033421,.073266){9}{\line(0,1){.073266}}
\multiput(29.693,27.803)(-.03336,.064342){10}{\line(0,1){.064342}}
\multiput(29.36,28.446)(-.033232,.056892){11}{\line(0,1){.056892}}
\multiput(28.994,29.072)(-.033049,.050552){12}{\line(0,1){.050552}}
\multiput(28.598,29.678)(-.0328165,.0450689){13}{\line(0,1){.0450689}}
\multiput(28.171,30.264)(-.0325401,.0402631){14}{\line(0,1){.0402631}}
\multiput(27.715,30.828)(-.0322235,.0360028){15}{\line(0,1){.0360028}}
\multiput(27.232,31.368)(-.03187,.0321894){16}{\line(0,1){.0321894}}
\multiput(26.722,31.883)(-.0356796,.032581){15}{\line(-1,0){.0356796}}
\multiput(26.187,32.372)(-.0399365,.03294){14}{\line(-1,0){.0399365}}
\multiput(25.628,32.833)(-.0447394,.0332644){13}{\line(-1,0){.0447394}}
\multiput(25.046,33.265)(-.05022,.033551){12}{\line(-1,0){.05022}}
\multiput(24.444,33.668)(-.051845,.030982){12}{\line(-1,0){.051845}}
\multiput(23.821,34.04)(-.058187,.030909){11}{\line(-1,0){.058187}}
\multiput(23.181,34.38)(-.065636,.030735){10}{\line(-1,0){.065636}}
\multiput(22.525,34.687)(-.074555,.030437){9}{\line(-1,0){.074555}}
\multiput(21.854,34.961)(-.085492,.029977){8}{\line(-1,0){.085492}}
\multiput(21.17,35.201)(-.099305,.029298){7}{\line(-1,0){.099305}}
\multiput(20.475,35.406)(-.117429,.028307){6}{\line(-1,0){.117429}}
\multiput(19.77,35.576)(-.17806,.03354){4}{\line(-1,0){.17806}}
\put(19.058,35.71){\line(-1,0){.7181}}
\put(18.34,35.808){\line(-1,0){.7221}}
\put(17.618,35.87){\line(-1,0){.7243}}
\put(16.894,35.895){\line(-1,0){.7247}}
\put(16.169,35.884){\line(-1,0){.7232}}
\put(15.446,35.837){\line(-1,0){.7199}}
\put(14.726,35.753){\line(-1,0){.7148}}
\multiput(14.011,35.633)(-.141565,-.031151){5}{\line(-1,0){.141565}}
\multiput(13.303,35.477)(-.116515,-.031864){6}{\line(-1,0){.116515}}
\multiput(12.604,35.286)(-.098369,-.032303){7}{\line(-1,0){.098369}}
\multiput(11.916,35.06)(-.084541,-.032561){8}{\line(-1,0){.084541}}
\multiput(11.239,34.8)(-.073596,-.032689){9}{\line(-1,0){.073596}}
\multiput(10.577,34.505)(-.064671,-.032716){10}{\line(-1,0){.064671}}
\multiput(9.93,34.178)(-.057221,-.032663){11}{\line(-1,0){.057221}}
\multiput(9.301,33.819)(-.050879,-.032543){12}{\line(-1,0){.050879}}
\multiput(8.69,33.428)(-.045394,-.0323654){13}{\line(-1,0){.045394}}
\multiput(8.1,33.008)(-.0405856,-.0321369){14}{\line(-1,0){.0405856}}
\multiput(7.532,32.558)(-.0363223,-.0318629){15}{\line(-1,0){.0363223}}
\multiput(6.987,32.08)(-.0346727,-.0336505){15}{\line(-1,0){.0346727}}
\multiput(6.467,31.575)(-.0329352,-.0353529){15}{\line(0,-1){.0353529}}
\multiput(5.973,31.045)(-.0333367,-.039606){14}{\line(0,-1){.039606}}
\multiput(5.506,30.49)(-.033709,-.0444054){13}{\line(0,-1){.0444054}}
\multiput(5.068,29.913)(-.0314314,-.0460456){13}{\line(0,-1){.0460456}}
\multiput(4.659,29.314)(-.031497,-.051533){12}{\line(0,-1){.051533}}
\multiput(4.281,28.696)(-.031488,-.057876){11}{\line(0,-1){.057876}}
\multiput(3.935,28.059)(-.031388,-.065326){10}{\line(0,-1){.065326}}
\multiput(3.621,27.406)(-.031179,-.074248){9}{\line(0,-1){.074248}}
\multiput(3.341,26.738)(-.030828,-.085188){8}{\line(0,-1){.085188}}
\multiput(3.094,26.056)(-.030287,-.099008){7}{\line(0,-1){.099008}}
\multiput(2.882,25.363)(-.029477,-.117141){6}{\line(0,-1){.117141}}
\multiput(2.705,24.661)(-.028253,-.142172){5}{\line(0,-1){.142172}}
\put(2.564,23.95){\line(0,-1){.7171}}
\put(2.458,23.233){\line(0,-1){.7215}}
\put(2.389,22.511){\line(0,-1){1.4488}}
\put(2.36,21.062){\line(0,-1){.7236}}
\put(2.401,20.339){\line(0,-1){.7207}}
\put(2.477,19.618){\line(0,-1){.7159}}
\multiput(2.59,18.902)(.029738,-.141868){5}{\line(0,-1){.141868}}
\multiput(2.739,18.193)(.0307,-.116827){6}{\line(0,-1){.116827}}
\multiput(2.923,17.492)(.031321,-.098686){7}{\line(0,-1){.098686}}
\multiput(3.142,16.801)(.031717,-.084861){8}{\line(0,-1){.084861}}
\multiput(3.396,16.122)(.031953,-.073918){9}{\line(0,-1){.073918}}
\multiput(3.684,15.457)(.03207,-.064995){10}{\line(0,-1){.064995}}
\multiput(4.004,14.807)(.032091,-.057544){11}{\line(0,-1){.057544}}
\multiput(4.357,14.174)(.032034,-.051201){12}{\line(0,-1){.051201}}
\multiput(4.742,13.56)(.0319111,-.0457145){13}{\line(0,-1){.0457145}}
\multiput(5.157,12.965)(.0317305,-.0409041){14}{\line(0,-1){.0409041}}
\multiput(5.601,12.393)(.031499,-.0366383){15}{\line(0,-1){.0366383}}
\multiput(6.073,11.843)(.033303,-.0350066){15}{\line(0,-1){.0350066}}
\multiput(6.573,11.318)(.0350226,-.0332862){15}{\line(1,0){.0350226}}
\multiput(7.098,10.819)(.0392716,-.0337301){14}{\line(1,0){.0392716}}
\multiput(7.648,10.346)(.0409194,-.0317109){14}{\line(1,0){.0409194}}
\multiput(8.221,9.902)(.0457299,-.0318891){13}{\line(1,0){.0457299}}
\multiput(8.815,9.488)(.051217,-.032009){12}{\line(1,0){.051217}}
\multiput(9.43,9.104)(.057559,-.032063){11}{\line(1,0){.057559}}
\multiput(10.063,8.751)(.06501,-.032038){10}{\line(1,0){.06501}}
\multiput(10.713,8.431)(.073933,-.031918){9}{\line(1,0){.073933}}
\multiput(11.379,8.143)(.084877,-.031676){8}{\line(1,0){.084877}}
\multiput(12.058,7.89)(.098701,-.031273){7}{\line(1,0){.098701}}
\multiput(12.748,7.671)(.116841,-.030644){6}{\line(1,0){.116841}}
\multiput(13.45,7.487)(.141883,-.02967){5}{\line(1,0){.141883}}
\put(14.159,7.339){\line(1,0){.716}}
\put(14.875,7.226){\line(1,0){.7207}}
\put(15.596,7.15){\line(1,0){.7237}}
\put(16.319,7.11){\line(1,0){.7248}}
\put(17.044,7.107){\line(1,0){.724}}
\put(17.768,7.14){\line(1,0){.7214}}
\put(18.489,7.209){\line(1,0){.717}}
\multiput(19.207,7.315)(.142158,.028321){5}{\line(1,0){.142158}}
\multiput(19.917,7.457)(.117127,.029533){6}{\line(1,0){.117127}}
\multiput(20.62,7.634)(.098993,.030335){7}{\line(1,0){.098993}}
\multiput(21.313,7.846)(.085174,.030869){8}{\line(1,0){.085174}}
\multiput(21.994,8.093)(.074233,.031214){9}{\line(1,0){.074233}}
\multiput(22.663,8.374)(.065311,.03142){10}{\line(1,0){.065311}}
\multiput(23.316,8.688)(.057861,.031515){11}{\line(1,0){.057861}}
\multiput(23.952,9.035)(.051518,.031522){12}{\line(1,0){.051518}}
\multiput(24.57,9.413)(.0460305,.0314535){13}{\line(1,0){.0460305}}
\multiput(25.169,9.822)(.0443892,.0337303){13}{\line(1,0){.0443892}}
\multiput(25.746,10.261)(.03959,.0333557){14}{\line(1,0){.03959}}
\multiput(26.3,10.727)(.035337,.0329522){15}{\line(1,0){.035337}}
\multiput(26.83,11.222)(.0336338,.0346889){15}{\line(0,1){.0346889}}
\multiput(27.335,11.742)(.0318454,.0363377){15}{\line(0,1){.0363377}}
\multiput(27.812,12.287)(.0321174,.0406011){14}{\line(0,1){.0406011}}
\multiput(28.262,12.856)(.0323436,.0454095){13}{\line(0,1){.0454095}}
\multiput(28.682,13.446)(.032519,.050895){12}{\line(0,1){.050895}}
\multiput(29.073,14.057)(.032636,.057237){11}{\line(0,1){.057237}}
\multiput(29.432,14.686)(.032685,.064687){10}{\line(0,1){.064687}}
\multiput(29.758,15.333)(.032653,.073611){9}{\line(0,1){.073611}}
\multiput(30.052,15.996)(.032521,.084557){8}{\line(0,1){.084557}}
\multiput(30.312,16.672)(.032256,.098384){7}{\line(0,1){.098384}}
\multiput(30.538,17.361)(.031808,.11653){6}{\line(0,1){.11653}}
\multiput(30.729,18.06)(.031083,.14158){5}{\line(0,1){.14158}}
\put(30.885,18.768){\line(0,1){.7148}}
\put(31.004,19.483){\line(0,1){.7199}}
\put(31.088,20.203){\line(0,1){1.2974}}
\put(16.75,36){\circle*{2}}
\put(7,31.75){\circle*{2}}
\put(2.75,25.25){\circle*{2}}
\put(3.25,16.75){\circle*{2}}
\put(7,11){\circle*{2}}
\put(13.75,8){\circle*{2}}
\put(22,8.75){\circle*{2}}
\put(27,12.25){\circle*{2}}
\put(30.5,21.25){\circle*{2}}
\put(26.5,32.25){\circle*{2}}
\qbezier(26.5,32.5)(22,53.25)(41.5,50)
\qbezier(41.5,50)(49.625,47.875)(52.25,33.25)
\qbezier(52.25,33.25)(57.125,7.5)(30.5,20.75)
\put(25.75,39){\circle*{2}}
\put(29.5,48.5){\circle*{2}}
\put(40.75,50.25){\circle*{2}}
\put(48.5,44.5){\circle*{2}}
\put(51.25,36.5){\circle*{2}}
\put(53,27.75){\circle*{2}}
\put(51,19){\circle*{2}}
\put(39,17.25){\circle*{2}}

\linethickness{1pt}\qbezier(30.75,20.75)(31.5,26)(26.25,32.25)\linethickness{0.4pt}
\multiput(27.25,28)(.6875,.03125){8}{\line(1,0){.6875}}
\multiput(32.75,28.25)(-.03125,.0625){8}{\line(0,1){.0625}}
\multiput(30,32)(-.0333333,-.2416667){30}{\line(0,-1){.2416667}}
\qbezier(73,27.25)(58.375,95.5)(118.25,94.75)
\qbezier(118.25,94.75)(194.875,93.875)(175,59.5)
\qbezier(175,59.5)(152.875,20.5)(130.25,11.5)
\qbezier(130.25,11.5)(80.25,-4.875)(73.25,27.25)
\qbezier(73.5,72)(84.875,41.25)(108.75,60.5)
\qbezier(108.75,60.5)(129.125,77.25)(101,93)
\multiput(76.25,83.75)(.17083333,-.03333333){60}{\line(1,0){.17083333}}
\multiput(81.25,88.75)(.03125,-1.3125){8}{\line(0,-1){1.3125}}
\end{picture}

\subsection{Embedding into 2-generated groups}

 \begin{theorem} (Sapir, \v Spakulov\' a \cite{SS})\label{SS1} Consider a group $G = \langle x_1,x_2,\dots,x_k | R = 1\rangle$, where $R$ is a word in the free group on $\{x_1,x_2,\dots,x_k\}$, $k\geq2$.
Assume the sum of exponents of $x_k$ in $R$ is zero and that the maximal Magnus $x_k$-index of $x_1$ is unique. Then $G$ can be embedded into a 2-generated 1-relator group which is an ascending HNN extension of a finitely generated free group.
\end{theorem}

The embedding is given by the map $x_i\mapsto w_i$, $i=1,...,k$ where

\begin{align*} w_1&=aba^2b...a^nba^{n+1}ba^{-n-1}ba^{-n}b...a^{-2}ba^{-1}b\\ w_i&=ab^ia^2b^i...a^nb^ia^{-n}b^i...a^{-2}b^ia^{-1}b^i, \quad \mathrm{ for}\quad 1<i<k\\ w_k&=ab^ka^2b^k...a^nb^ka^{-n}b^k...a^{-2}b^k \end{align*}
The injectivity of that map follows from Theorem \ref{O1}.

\section{Probability Theory. Brownian motions}

Let $C$ be the space of all continuous functions $f\colon [0,+\infty]\to \R^k$ with $f(0)=0$. We can define a $\sigma$-algebra structure on that space generated by the sets of functions of the form $U(t_1,x_1, t_2, x_2, ..., t_n, x_n)$ where $t_i\in [0,+\infty], x_i\in \R^k$. This set consists of all functions $f\in C$ such that $f(t_i)=x_i$. A measure $\mu$ on $C$ is called the {\em Wiener's measure} if for every Borel set $A$ in $\R^k$ and every $t< s \in [0,+\infty]$ the probability that $f(t)-f(s)$ is in $A$ is $$\frac{1}{\sqrt{2\pi(t-s)}}\int_A e^\frac{-|x|^2}{2(t-s)} dx.$$ This means that  (by definition) Brownian motion is a continuous Markov stationary process with normally distributed increments.

\subsection{Donsker's theorem (modified)}

A standard tool in dealing with random walks is to consider rescaled limits of them. The rescaled limit of random walks is ``usually" a Brownian motion. This is the case in our situation too. Cyclically reduced relators correspond to cyclically reduced walks, i.e. walks without backtracking such that the labels of the first and the last steps are not mutually inverse.

\begin{theorem}(Sapir, \v Spakulov\' a, \cite{SS}) Let $P^{CR}_n$ be the uniform distribution on the set of cyclically reduced random walks of length $n$ in $\R^k$.  Consider a piecewise linear function $Y_n(t):[0,1]\to \R^k$, where the line segments are connecting points $Y_n(t)=S_{nt}/\sqrt{n}$ for $t=0,1/n,2/n$, $\dots,n/n=1$, where $(S_n)$ has a distribution according to $P^{CR}_n$. Then $Y_n(t)$ converges in distribution to a Brownian motion, as $n \to \infty$.
\end{theorem}

The main ingredient in the proof is Rivin's Central Limit Theorem for cyclically reduced walks \cite{Rivin}.

\subsection{Convex hull of Brownian motion and maximal Magnus indices.}

Let again $w$ be the walk in $\Z^k$ corresponding to the relator $R$ in $k$ generators. Suppose that it connects $O$ and $M$. Consider the hyperplane $P$ that is orthogonal to $\vec{OM}$, the projection $w'$ of $w$ onto $P$, and the convex hull $\Delta$ of that projection.  From Theorem \ref{SS1}, it follows that the 1-relator group $G$ is inside an ascending HNN extension of a free group if there exists a vertex of $\Delta$ that is visited only once by $w'$. The idea to prove that this happens with probability tending to 1 is the following.

{\bf Step 1.} We prove that the number of vertices of $\Delta$ is growing (a.s.) with the length of $w$ (here we use that $k\ge 3$; if $k=2$, then $\Delta$ has just two vertices). Indeed, if the number of vertices is bounded with positive probability, then with positive probability the limit of random walks $w'$ (which is a Brownian bridge) would have non-smooth convex hull which is impossible by a theorem about Brownian motions (Theorem of Cranston-Hsu-March \cite{CHM}, 1989).

{\em Step 2.}  For every vertex of $\Delta$ for any  `bad" walk $w'$ or length $r$ we construct (in a bijective manner) a ``good" walk $w'$ of length $r+4$. This implies that the number of vertices of ``bad" walks is bounded almost surely if the probability of a ``bad" walk is $>0$.

Here is the walk in $\Z^3$ corresponding to the word
$$cb^{-1}acac^{-1}b^{-1}caca^{-1}b^{-1}aab^{-1}c.$$
and its projection onto the plane perpendicular to $\vec{OM}$. The second (left bottom) picture shows an approximal projection. One can see that the walk is bad: every vertex of the convex hull of the projection is visited twice.
\begin{center}
\includegraphics[width=200pt]{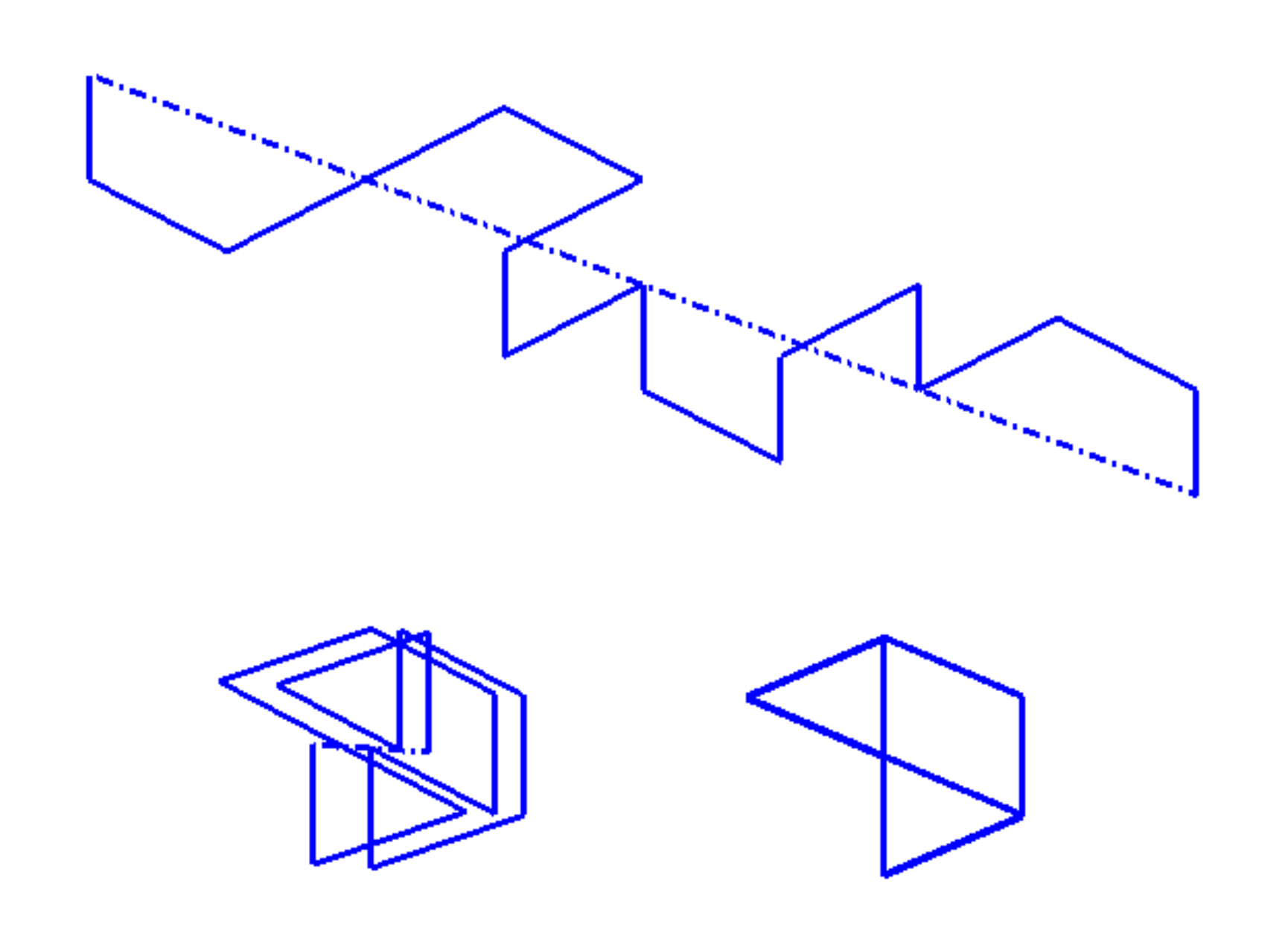}
\end{center}

The walk corresponds to 5 different good walks (5 is the number of vertices in the convex hull) by inserting squares after the second visit of the vertex of the convex hull. Here is the walk and its projection corresponding to the word $$cb^{-1}acac^{-1}b^{-1}caca^{-1}b^{-1}
((b^{-1}cbc^{-1}))aab^{-1}c.$$
\begin{center}
\includegraphics[width=200pt]{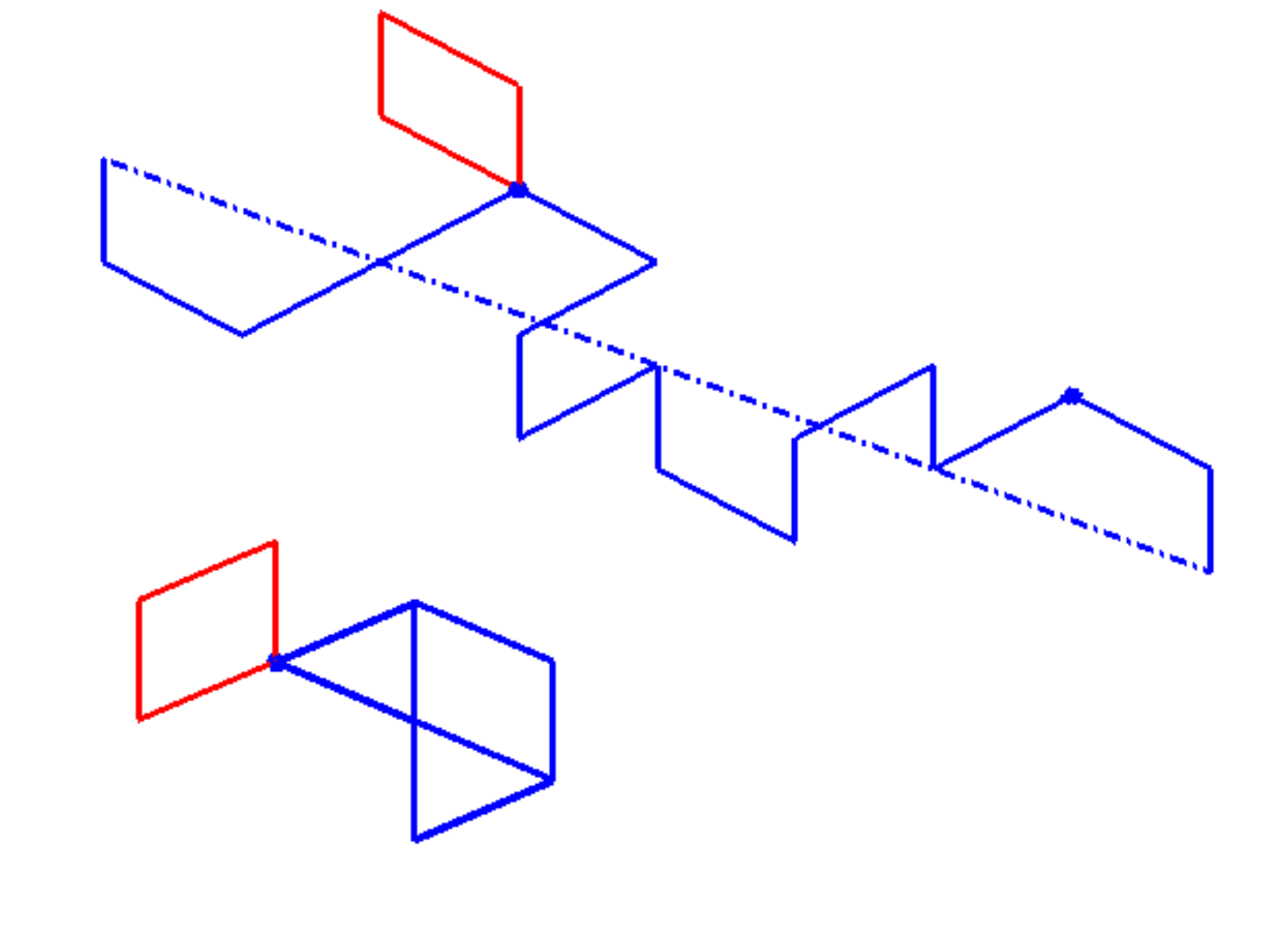}
\end{center}

Thus almost surely every 1-relator group with at least 3 generators is inside a 2-generated 1-relator group which is an ascending HNN extension of a free group. This implies, by Feighn and Handel \cite{FH}, that a.s. every 1-relator group with at least 3 generators is coherent. Residual properties need new ideas.

\section{Algebraic geometry}

\subsection{Periodic points of a word map}

Consider the group $$G=\la x_1,...,x_k,t\mid x_1^t=\phi(x_1), ..., x_k^t=\phi(x_k)\ra$$ for some injective endomorphism $\phi$ of $F_k$.

For example, $$H_T=\la x,y,t\mid txt\iv =xy, tyt\iv=yx\ra.$$ So the endomorphism $$\phi\colon x\mapsto xy, y\mapsto yx.$$
It is easy to see that every element in $G$ has the form $t^kw(x_1,...,x_k)t^l$. If $k+l\ne 0$, then the element survives the homomorphism $x_i\mapsto 0$, $t\mapsto 1$ $G\to \Z$. Hence we can assume that $k+l=0$, and the element is conjugate to $w(x_1,...,x_k)$. Thus it is enough to consider elements from $F_k$ only.  Consider any $w=w(x,y)\ne 1$.
We want to find $\psi\colon G\to V$ with $\psi(w)\ne 1$, $|V|<\infty$.
Suppose that $\psi$ exists.

Let us denote $\psi(x), \psi(y), \psi(t)$ by $\bar x,\bar y,\bar t$. So we want $$w(\bar x, \bar y)\ne 1.$$ Note: $$\bar t(\bar x,\bar y)\bar t\iv = (\bar x\bar y,\bar y\bar x)=(\phi(x), \phi(y))$$ We can continue: $$\bar t^2(\bar x,\bar y) \bar t^{-2} = (\bar x\bar y\bar y\bar x,\bar y\bar x\bar x\bar y)= (\phi^2(\bar x), \phi^2(\bar y)).$$ $$...$$ $$\bar t^k(\bar x,\bar y)\bar t^{-k} =(\phi^k(\bar x), \phi^k(\bar y)).$$

Since $\bar t$ has finite order in $V$, for some $k$, we must have $$(\phi^k(\bar x),\phi^k(\bar y))=(\bar x,\bar y).$$ So $(\bar x, \bar y)$ is a periodic point of the map $$\tilde\phi\colon (a,b)\mapsto (ab,ba).$$ on the ``space" $V\times V$.

So if $G$ is residually finite then for every $w(x,y)\ne 1$, we found a finite group $V$ and a periodic point $(\bar x, \bar y)$ of the map $$\tilde\phi\colon (a,b)\mapsto (\phi(a), \phi(b))$$ on $V\times V$ such that $$w(\bar x, \bar y)\ne 1.$$ So the periodic point should be outside the ``subvariety" given by $w=1$.

{\bf The key observation.} If $(\bar x, \bar y)$ is periodic with period of length $\ell$, $w(\bar x, \bar y)\ne 1$, then there exists a homomorphism from $G$ into the wreath product $V'=V \wr \Z/\ell\Z=(V\times V \times ... \times V) \leftthreetimes \Z/\ell\Z$ which separates $w$ from 1.

Indeed, the map $$x\mapsto ((\bar x, \phi(\bar x), \phi^2(\bar x), ..., \phi^{l-1}(\bar x)),0),$$ $$y\mapsto ((\bar y, \phi(\bar y), \phi^2(\bar y), ..., \phi^{\ell-1}(\bar y)),0),$$ $$t\mapsto ((1,1,...,1),{\bf 1})$$
extends to a homomorphism $\gamma:G\to V'$
and $$\gamma(w)=((w(\bar x, \bar y),...), 0)\ne 1.$$

{\bf The idea} Thus in order to prove that the group $\HNN_\phi(F_k)$ is residually finite, we need, for every word $w\ne 1$ in $F_k$, find a finite group $G$ and a periodic point of the map $\tilde\phi\colon G^k\to G^k$ outside the ``subvariety" given by the equation $w=1$.

Consider again the group $H_T=\la a,b, t\mid tat\iv=ab, tbt\iv=ba\ra$. Consider two matrices $$U=\left[\begin{array}{ll}1& 2\\ 0& 1\end{array}\right], V=\left[\begin{array}{ll}1& 0\\ 2& 1\end{array}\right].$$ They generate a free subgroup in $\SL_2(\mathbb{Z})$ (Sanov \cite{Sanov}). Then the matrices $$A=UV=\left[\begin{array}{ll}5& 2\\ 2& 1\end{array}\right], B=VU=\left[\begin{array}{ll}1& 2\\ 2& 5\end{array}\right]$$ also generate a free subgroup. Now let us iterate the map $\psi\colon (x,y)\to (xy,yx)$ starting with $(A,B)$ mod 5 (that is we are considering the finite group $\SL_2(\Z/5\Z)$):

$$\left(\left [\begin {array}{cc} 5&2\\\noalign{\medskip}2&1\end {array} \right ], \left [\begin {array}{cc} 1&2\\\noalign{\medskip}2&5\end {array} \right ]\right)\to \left(\left [\begin {array}{cc} 4&0\\\noalign{\medskip}4&4\end {array} \right ], \left [\begin {array}{cc} 4&4\\\noalign{\medskip}0&4\end {array} \right ]\right)\to$$

$$\left(\left [\begin {array}{cc} 5&2\\\noalign{\medskip}2&1\end {array} \right ], \left [\begin {array}{cc} 1&2\\\noalign{\medskip}2&5\end {array} \right ]\right)\to \left(\left [\begin {array}{cc} 4&0\\\noalign{\medskip}4&4\end {array} \right ], \left [\begin {array}{cc} 4&4\\\noalign{\medskip}0&4\end {array} \right ]\right)\to$$

$$\left(\left [\begin {array}{cc} 1&1\\\noalign{\medskip}1&2\end {array} \right ], \left [\begin {array}{cc} 2&1\\\noalign{\medskip}1&1\end {array} \right ]\right)\to \left(\left [\begin {array}{cc} 3&2\\\noalign{\medskip}4&3\end {array} \right ], \left [\begin {array}{cc} 3&4\\\noalign{\medskip}2&3\end {array} \right ]\right)\to$$

$$\left(\left [\begin {array}{cc} 3&3\\\noalign{\medskip}3&0\end {array} \right ], \left [\begin {array}{cc} 0&3\\\noalign{\medskip}3&3\end {array} \right ]\right)\to \left(\left [\begin {array}{cc} 4&3\\\noalign{\medskip}0&4\end {array} \right ], \left [\begin {array}{cc} 4&0\\\noalign{\medskip}3&4\end {array} \right ]\right)\to$$

$$ \left(\left [\begin {array}{cc} 5&2\\\noalign{\medskip}2&1\end {array} \right ], \left [\begin {array}{cc} 1&2\\\noalign{\medskip}2&5\end {array} \right ]\right).$$

Thus the point $(A,B)$ is periodic in $\SL_2(\Z/5\Z)$ with period 6.

\subsection{Dynamics of polynomial maps over local fields I}

Let us replace $5$ by $25$, $125$, etc. It turned out that $(A,B)$ is periodic in $SL_2(\Z/25\Z)$ with period 30, in $\SL_2(\Z/125\Z)$ with period $150$, etc.

Moreover there exists the following general Hensel-like statement:

\begin{theorem}\label{Hensel} Let $P\colon \Z^n\to \Z^n$ be a polynomial map with integer coefficients. Suppose that a point $\vec x$ is periodic with period $d$ modulo some prime $p$, and the Jacobian $J_P(x)$ is not zero. Then $\vec x$ is periodic modulo $p^k$ with period $p^{k-1}d$ for every $k$.
\end{theorem}

This theorem does not apply to our situation straight away because it is easy to see that the Jacobian of our map is 0, but slightly modifying the map (decreasing the dimension), we obtain the result.

Now take any word $w\ne 1$ in $x,y$. Since $\la A,B\ra$ is free in $\SL_2(\Z)$, the matrix $w(A,B)$ is not 1, and there exists $k\ge 1$ such that $w(A,B)\ne 1 \mod 5^k$.

Therefore our group $H_T$ is residually finite.

\subsection{Reduction to polynomial maps over finite fields}

Let us try to generalize the example of $H_T$.

Consider an arbitrary ascending HNN extension of a free group  $G=\la a_1,...,a_k,t \mid ta_it\iv=w_i, i=1,...,k\ra$. Consider the ring of matrices $M_2(\Z)$.

The map $\psi\colon M_2(\Z)^k\to M_2(\Z)^k$ is given by $$\vec x\mapsto (w_1(\vec x), ...,w_k(\vec x)).$$

It can be considered as a
polynomial map $A^{4k}\to A^{4k}$ (where $A^l$ is the affine space of dimension $l$). In fact there is a slight problem here with the possibility that $w_i$ may contain inverses of elements. To resolve it,  we replace inverses by the adjoint operation applicable to all matrices. For invertible matrices the inverse and the adjoint differ by a scalar multiple. Thus the group we consider will be $\mathrm{PGL}_2(.)$ instead of $\mathrm{SL}_2(.)$. Here "." is any finite field of, say, characteristic $p$.

Thus our problem is reduced to the following:

{\bf Problem.} Let $P$ be a polynomial map $A^n\to A^n$ with integer coefficients. Show that the set of periodic points of $P$ is Zariski dense.

\subsection{Deligne problem}

Fixed points are not enough.

\medskip

{\bf Example:} $x\mapsto x+1$ does not have fixed points.

Deligne suggested to consider consider quasi-fixed points,
$$P(x)=x^{p^n} (=Fr^n(x))$$ where $Fr$ is the Frobenius map (raising all coordinates to the power $p$).

Note that all quasi-fixed points are periodic because $Fr$ is the automorphism of finite order, and commutes with $P$ since all coefficients of $P$ are integers.

{\bf Deligne conjecture}  proved by Fujiwara and Pink: If $P$ is dominant and quasi-finite then the set of quasi-fixed points is Zariski dense. Unfortunately in our case the map is rarely dominant or quasi-finite. Say, in the case of $H_T$, the matrices $xy, yx$ have the same trace, so they satisfy a polynomial equation that does not hold for the pair $(x,y)$. Hence the map is not dominant.

\subsection{The main results}

\begin{theorem} \label{BS2} (Borisov, Sapir \cite{BS1}) Let $P^{n}\colon A^n(\F_q)\to A^n(\F_q)$ be the $n$-th
iteration of $P$. Let $V$ be the Zariski closure of $P^{n}(A^n)$. The set of its geometric points is $V(\overline{\F_q}),$ where $\overline{\F_q}$ is the algebraic closure of $\F_q.$ Then the following holds. \begin{enumerate} \item All quasi-fixed points of $P$ belong to $V(\overline{\F_q}).$
\item Quasi-fixed points of $P$ are Zariski dense in $V$. In other words, suppose $W\subset V$ is a proper Zariski closed subvariety of $V$. Then for some $Q=q^m$ there is a point $$(a_1,...,a_n)\in V(\overline{\F_q})\setminus W(\overline{\F_q})$$ such that 
    
\begin{equation}\label{eq3} \left\{ \begin{array}{l} f_1(a_1,..., a_n)=a_1^Q\\ f_2(a_1,...,a_n)=a_2^Q\\ ...\\ f_n(a_1,...,a_n)=a_n^Q
 \end{array}\right.\end{equation}
\end{enumerate}
\end{theorem}

This, as before, implies

\begin{theorem}\label{BS3} (Borisov, Sapir \cite{BS1}). Every ascending HNN extension of a free group is residually finite.
\end{theorem}

The proof of Theorem \ref{BS2} from \cite{BS1} is non-trivial but relatively short. Unfortunately the naive approach based on the B\'ezout theorem fails. B\'ezout theorem (that the number of solutions of the system of equations (\ref{eq3}) is $Q^n$ for large enough $Q$ \cite{Fu}) gives only the number of solutions of our system of equations but does not tell us that there are solutions outside a given subvariety (in fact all solutions can, in principle, coincide, as for the equation $x^Q=0$). Our proof gives a very large $Q$ and uses manipulations with ideals in the ring of polynomials, and its localizations.

Here is a more detailed description of the proof of Theorem \ref{BS3}.

We denote by $I_Q$ the ideal in $\overline{\F_q}[x_1,...,x_n]$ generated by the polynomials $f_i(x_1,...,x_n)-x_i^Q,$ for $i=1,2,...,n.$

{\bf Step 1.} For a big enough $Q$ the ideal $I_Q$ has finite codimension in the ring $\overline{\F_q}[x_1,..x_n]$.

{\bf Step 2.}  For all $1\leq i\leq n$ and $j\geq 1$ $$f_i^{(j)}(x_1,...,x_n)-x_i^{Q^j}\in I_Q.$$

{\bf Step 3.} There exists a number $k$ such that for every quasi-fixed point $(a_1,...,a_n)$ with big enough $Q$ and for every $1\leq i\leq n$ the polynomial $$(f_i^{(n)}(x_1,...,x_n)- f_i^{(n)}(a_1,...,a_n))^k$$ is contained in the localization of $I_Q$ at $(a_1,...,a_n).$

Let us fix some polynomial $D$ with the coefficients in a finite extension of $\F_q$ such that it vanishes on $W$ but not on $V$.

{\bf Step 4. } There exists a positive integer $K$ such that for all quasi-fixed points $(a_1,...,a_n)\in W$ with big enough $Q$ we get $$R=(D(f_1^{(n)}(x_1,...,x_n),...,f_n^{(n)}(x_1,...,x_n)))^K\equiv 0 (\mod I_Q^{(a_1,...,a_n)})$$

We know that all points with $P(x)=x^Q$ belong to $V$. We want to prove that some of them do not belong to $W$. We suppose that they all do, and we are going to derive a contradiction.

{\bf Step 5.} First of all, we claim that in this case $R$ lies in the localizations of $I_Q$ with respect to all maximal ideals of the ring of polynomials.

This implies that $R\in I_Q$.
Therefore there exist polynomials $u_1,...u_n$ such that \begin{equation}\label{eq90} R=\sum \limits_{i=1}^{n} u_i \cdot (f_i-x_i^Q) \end{equation}

{\bf Step 6.} We get a set of $u_i$'s with the following property: \begin{itemize} \item For every $i<j$ the degree of $x_i$ in every monomial in $u_j$ is smaller than $Q$. \end{itemize}

{\bf Step 7.} We look how the monomials cancel in the equation (1) and get a contradiction.

\subsection{Hrushovsky's result}

Hrushovsky \cite{H} managed to replace $A^n$ in our statement by arbitrary variety $V$. His proof is very non-rivial and uses model theory as well as deep algebraic geometry. Using his result we prove

\begin{theorem}\label{BS4} (Borisov, Sapir, \cite{BS1}) The ascending HNN extension of any finitely generated linear group is residually finite.
\end{theorem}

This theorem is stronger than Theorem \ref{BS3} because free groups are linear. It also applies to, say, right angled Artin groups (which usually have many injective endomorphisms).

For non-linear residually finite groups the statement of Theorem \ref{BS4} is not true: the Lysenok HNN extension of Grigorchuk's group \cite{Ly} is not residually finite  (see Sapir and Wise \cite{SW}).

\subsection{Dynamics of polynomial maps over local fields II}

It remains to show that the ascending HNN extensions of free groups are virtually residually (finite $p$-groups) for almost all primes $p$.

{\bf Example} Consider again the group $H_T=\la a,b, t\mid tat\iv=ab, tbt\iv=ba\ra$.
Let us prove that $G$ is virtually residually ($5$-group).
We know that the pair of  matrices $A=\left[\begin{array}{ll}5& 2\\
2& 1\end{array}\right], B=\left[\begin{array}{ll}1& 2\\ 2& 5\end{array}\right]$ generates a free subgroup and is a periodic point of the map $\psi\colon (x,y)\to (xy,yx)$ modulo $5^d$. The period is
 $\ell_d=6*5^{d-1}$ by Theorem \ref{Hensel}.

Then the group $H_T$ is approximated by subgroups of the finite groups $SL_2(\Z/5^d\Z)^{\ell_d}\leftthreetimes \Z/\ell_d\Z$. Let $\nu_d$ be the corresponding homomorphisms.

 Let $G_d=SL_2(\Z/5^d\Z)$. There exists a natural homomorphism  $\mu_d\colon G_d^{\ell_d} \to G_1^{\ell_d}$. The kernel is a $5$-group.

The image $\langle x_d, y_d\rangle$ of $F_2=\langle
x,y\rangle$ under $\nu_d\mu_d$ is inside the direct power
$G_1^{\ell_d}$, hence $\nu_d\mu_d(F_2)$ is a 2-generated group in the
variety of groups generated by $G_1$, hence it is of bounded size, say, $M$.

Then there exists a  characteristic subgroup of $\mu_d(F_2)$ of index $\le M_1$ which is residually $5$-group.

The image of $H_T$ under $\mu_d$ is is an extension of a subgroup that has a $5$-subgroup of index $\le M_1$ by a cyclic group which is an image of $\Z/6*5^{d-1}\Z$.

Since $\Z/6*5^{d-1}\Z$ has
a 5-subgroup of index 6,
$\mu_d\nu_d(H)$ has a 5-subgroup of index
at most some constant $M_2$ (independent of $d$). Hence $H_T$ has a subgroup of index at most $M_2$
which is residually (finite $5$)-group. Thus $H_d$ is virtually residually (finite $5$-group).

Consider now the general situation.

To show that $G=\HNN_\phi(F_k)$ is virtually residually (finite $p$)-group for almost all $p$
we do the following:

\begin{itemize}
\item Instead of a tuple of matrices $(A_1,A_2,...,A_k)$ from $\SL_2(\Z)$ such that $(A_1 \mod p,..., A_k \mod
p)$ is periodic for the map $\phi$, we find a periodic tuple of matrices in $\SL_2(\OO/p\OO)$ where $\OO$ is the ring of integers of some finite extension of $\Q$ unramified at $p$ (this is possible to do using Theorem \ref{BS2}).

\item Then using a version of Theorem \ref{Hensel}, and a result of Breuillard and  Gelander \cite{BG} we lift
the matrices $A_1,..., A_n$ to the $p$-adic completion of $SL_2(\OO)$ in such a way that the lifts generate a free subgroup.
\end{itemize}

This completes the proof of Theorem \ref{thmain}.

\section{Applications}

\subsection{An application to pro-finite groups}

Theorem \ref{BS2} is equivalent to the following statement about pro-finite completions of free groups.

\begin{theorem}\label{BS6}(Borisov, Sapir \cite{BS1}) For every injective homomorphism $\phi\colon F_k\to F_k$ there exists a pro-finite completion $F$ of $F_k$, and an automorphism $\phi'$ of $F$ such that $\phi'_{F_k}=\phi$.
\end{theorem}

\proof Consider the homomorphisms $\mu_n$ from $\HNN_\phi(F_k)$ to finite groups.
The images $\mu_n(t)$ ($t$ is the free generator of $\HNN_\phi(F_k)$) induce automorphisms of the images of $F_k$. The intersections of kernels of $\mu_n$ with $F_k$ give a sequence of subgroups of finite index.
The corresponding pro-finite completion of $F_k$ is what we need. \endproof

\subsection{Some strange linear groups} \label{str}

Let $\HNN_\phi(F_k)=\la F_k, t\mid F_k^t=\phi(F_k)\ra$. Consider the normal closure $N$ of $F_k$ in $\HNN_\phi(F_k)$. Since $F_k^t\subseteq F_k$, we have $F_k\subseteq F_k^{t^{-1}}$:

\unitlength .6 mm 
\linethickness{0.4pt}
\ifx\plotpoint\undefined\newsavebox{\plotpoint}\fi 
\begin{picture}(204.25,92)(0,0)
\qbezier(67,40.5)(38.25,40.625)(67.5,14.25)
\qbezier(67.5,14.25)(94.875,38.875)(66.75,40)
\qbezier(66,59)(-9.25,57.75)(67.5,14.5)
\qbezier(67.5,14.5)(143.875,58.875)(65.75,58.75)
\put(66,62){\circle*{1.5}}
\put(66,66.75){\circle*{1.5}}
\put(66,71.5){\circle*{1.5}}
\qbezier(65,88.75)(-63.25,85.875)(67.5,14.5)
\qbezier(67.5,14.25)(204.25,92)(65,88.75)
\put(67,27.75){\makebox(0,0)[cc]{$F_k$}}
\put(66,49){\makebox(0,0)[cc]{$F_k^{t^{-1}}$}}
\put(92.25,73.5){\makebox(0,0)[cc]{$N=\cup_n F_k^{t^{-n}}$}}
\end{picture}

Hence the group $N$ is a locally free group where every finitely generated subgroup is inside a free group of rank $k$. If $\phi(F_k)\subset [F_k, F_k]$, then $[N,N]=N$. Nevertheless, $N$ is linear and in fact is inside $\SL_2(\C)$. This immediately follows from the proof of Theorem \ref{BS3}.

\section{Two possible approaches to Problem \ref{p1}. Some other open problems}

It was proved by Minasyan (unpublished) that the group $H_T$ is hyperbolic. Hence by \cite{Gr1, O1}, it would remain hyperbolic if we impose additional relation of the form $t^p=1$, where $t$ is the free generator of $H_T$, $p>>1$, a prime.

\begin{prob}\label{p8} Is the group $\la a,b, t\mid a^t=ab, b^t=ba, t^p=1\ra$ residually finite.
\end{prob}

If the answer is ``yes", then there are arbitrary large finite 2-generated groups $\la a,b\ra$ for which the pair $(a,b)$ is periodic with respect to the map $(x,y)\mapsto (xy,yx)$ with the fixed period $p$ (note that in the proof above the periods of points rapidly increase with $p$). We can add more torsion by imposing $a^p=1, b^p=1$ as well (but keeping the factor-group hyperbolic).

Let us also recall a problem by Olshanskii from Kourovka Notebook \cite[Problem 12.64]{KN} which can be slightly modified as follows.

 \begin{prob}
 Let $d\gg p\gg 1$ ($p$ a prime). Are there infinitely many 2-generated finite simple groups $G$ with generators $a,b$ such that for every word $w(x,y)$ of length at most $d$, $G$ satisfies $w(a,b)^p=1$? If the answer is ``no", then there exists a non-residually finite hyperbolic group. 
\end{prob}

Note that since (assuming the Classification of Finite Simple Groups) there are only 18 types of non-abelian finite simple groups, one can ask this question for each type separately. In particular, the specification to $A_n$ gives a nice, but still unsolved, combinatorial problem about permutations:
\begin{prob}
Are there two permutations $a$ and $b$, of $\{1,...n\}$ which together act transitively on the set, and such that $w(a,b)^p=1$ for every word $w$ of length $\le d$. Here $p$ and $d$ are sufficiently large (say, $d = \exp(\exp(p)), p \gg 1$), $n$ runs over an infinite set of numbers (the set depends on $p,d$). To exclude the case when $\la a, b\ra$ is a $p$-group (in which case we run into the restricted Burnside problem) we can assume that $p\not\mid n$.
\end{prob}

Another approach to Problem \ref{p1} is to consider multiple HNN extensions of free groups. Consider, for example the group $H_T(u,v)=\la a,b,t,s\mid a^t=ab, b^t=ba, a^s=u, b^s=v\ra$ where $u,v$ are words in $a,b$ ``sufficiently independent" from $ab,ba$. This group is again hyperbolic by, say, \cite{OlSq}.  Suppose that $H_T(u,v)$ has a homomorphism onto a finite group $G$, and $\bar a, \bar b$ are images of $a,b$, we assume $\bar a, \bar b\ne 1$. Then the pair $(\bar a, \bar b)$ is periodic for both maps $\phi: (x,y)\mapsto (xy,yx)$ and $\psi: (x,y)\mapsto (u(x,y), v(x,y))$. Moreover for every word $w(x,y)$, the pair $w(\phi, \psi)(a,b)$ is periodic both for $\phi$ and for $\psi$. If $\phi$ and $\psi$ are sufficiently independent, then it seems unlikely that these two maps have so many common periodic points. MAGMA  refuses to give any finite simple homomorphic image with $\bar a\ne 1, \bar b\ne 1$ for almost any choice of $u,v$.

Here is an ``easier" problem.

\begin{prob}\label{p5} Is the group $H_T=\la a,b,t\mid tat\iv=ab, tbt\iv=ba\ra$ linear?
\end{prob}

 Recall that the normal closure $N$ of $\{a,b\}$ is linear (see  Section \ref{str}). So $H_T$ is an extension of a linear group by a cyclic group. Still we conjecture that the answer to Problem \ref{p5} is ``No". One way to prove it would be to consider the minimal indices of subgroups of $H_T$ that are residually (finite $p$-groups). In linear groups these indices grow polynomially in terms of $p$. Our proof gives a much faster growing function (double exponential). Non-linear hyperbolic groups are known \cite{Kap1} but their presentations are much more complicated.

In \cite{DS}, we asked whether an ascending HNN extension $\HNN_\phi(F_k)$ where $\phi$ is a proper endomorphism, $k\ge 2$ can be inside $\SL_2(\C)$. Calegary and Dunfield constructed such an example \cite{CD}. But in their example, $\phi$ is reducible. Its components are an automorphism and a Baumslag-Solitar endomorphism of $F_1$ ($x\mapsto x^\ell$). Both Baumslag-Solitar groups and many finitely generated free-by-cyclic groups are inside $\SL_2(\C)$, so the result of \cite{CD} is not very surprising. Thus we still have the following

\begin{prob}\label{p6} Is there an irreducible proper injective endomorphism of $F_k$, $k\ge 2$, such that $\HNN_\phi(F_k)$ is embeddable into $\SL_2(\C)$?
\end{prob}

\begin{minipage}[t]{3 in}
\noindent Mark V. Sapir\\ Department of Mathematics\\
Vanderbilt University\\
m.sapir@vanderbilt.edu\\
http://www.math.vanderbilt.edu/$\sim$msapir\\
\end{minipage}


\begin{thebibliography}{kkkkl}

\bibitem[Ad]{NA} S. I. Adian, The Burnside problem and identities in
groups.  Ergebnisse der Mathematik und ihrer Grenzgebiete [Results
in Mathematics and Related Areas], {\bf 95}, Springer-Verlag,
Berlin-New York, 1979.

\bibitem[BP]{BP} Benjamin Baumslag, Stephen J. Pride, 
Groups with two more generators than relators.
J. London Math. Soc. (2) 17 (1978), no. 3, 425--426. 

\bibitem[BS1]{BS1} Alexander Borisov, Mark Sapir, Polynomial maps over finite fields
and residual finiteness of mapping tori of group endomorphisms,
Invent. Math., 160, 2 (2005), 341 - 356.


\bibitem[BS2]{BS2} Alexander Borisov, Mark Sapir, Polynomial maps over finite and $p$-adic fields and residual properties of ascending HNN extensions of free groups, International Mathematics Research Notices, 2009, 16, 3002-3015.


\bibitem[BG]{BG}Emanuel Breuillard, Tsachik Gelander,  A topological Tits alternative.  Ann. of Math. (2)  166  (2007),  no. 2, 427--474.

\bibitem[Br]{Brown}
Kenneth~S. Brown.
\newblock Trees, valuations, and the {B}ieri-{N}eumann-{S}trebel invariant.
\newblock {Invent. Math.}, 90(3):479--504, 1987.
\bibitem[DS]{DS} Cornelia Dru\c tu, Mark Sapir, Non-linear residually finite groups,
J. of Algebra, v. 284, 1 (2005), 174-178.

\bibitem[BM]{BM} Marc Burger, Shahar Mozes, Lattices in product of trees.
 Publ. IHES,  No. 92  (2000), 151--194 (2001).
		
\bibitem[CD]{CD} Danny Calegari, Nathan Dunfield,  An ascending HNN extension of a free group inside SL(2,C). Proc. Amer. Math. Soc. 134 (2006), no. 11, 3131-3136.

\bibitem[CHM]{CHM}
M.~Cranston, P.~Hsu, and P.~March.
\newblock Smoothness of the convex hull of planar {B}rownian motion.
\newblock {\em Ann. Probab.}, 17(1):144--150, 1989.

\bibitem[DT]{DT} Nathan M. Dunfield, Dylan P. Thurston,  A random tunnel number
one 3-manifold does not fiber over the circle. Geom. Topol. 10
(2006), 2431--2499.

\bibitem[FH]{FH} M. Feighn, M. Handel, Mapping tori of free group automorphisms
are coherent. Ann. Math. (2) 149, 1061--1077 (1999).

\bibitem[F]{Fu} William Fulton,  Intersection theory. Second edition. Ergebnisse
der Mathematik und ihrer Grenzgebiete. 3. Folge. A Series of
Modern Surveys in Mathematics [Results in Mathematics and Related
Areas. 3rd Series. A Series of Modern Surveys in Mathematics], 2.
Springer-Verlag, Berlin, 1998.

\bibitem[GMSW]{GMSW} Ross Geoghegan, Michael L. Mihalik,
Mark Sapir, Daniel T. Wise. Ascending HNN extensions of finitely
generated free groups are Hopfian. Bull. London Math. Soc. 33
(2001), no. 3, 292--298.

\bibitem[GNS]{GNS} R. I. Grigorchuk, V. V. Nekrashevich, V. I. Sushchanskii,  Automata, dynamical systems, and groups. Tr. Mat. Inst. Steklova  231  (2000),  Din. Sist., Avtom. i Beskon. Gruppy, 134--214.

\bibitem[Gr1]{Gr1} M. Gromov, Hyperbolic groups, in ''Essays in Group Theory'', S. M. Gersten, Ed..,
M.S.R.I. Pub. 8, pp. 75]263, Springer-Verlag, BerlinrNew York, 1987.

\bibitem[Hr]{H} E. Hrushovski. The Elementary Theory of the Frobenius
Automorphisms, arXive math.LO/0406514.




\bibitem[KSS]{KSS}
I. Kapovich, P. Schupp, V. Shpilrain,
\newblock Generic properties of {W}hitehead's algorithm and isomorphism
rigidity of random one-relator groups.
\newblock {\em Pacific J. Math.}, 223(1):113--140, 2006.

\bibitem[M.Kap]{Kap1} Michael Kapovich, Representations of polygons of finite groups.  Geom. Topol.  9  (2005), 1915--1951.

\bibitem[KN]{KN} The Kourovka notebook. Unsolved problems in group theory. Sixteenth edition. Including archive of solved problems. Edited by V. D. Mazurov and E. I. Khukhro. Russian Academy of Sciences Siberian Division, Institute of Mathematics, Novosibirsk, 2006. 178 pp.

\bibitem[La]{La} Marc Lackenby,
Adding high powered relations to large groups.
Math. Res. Lett. 14 (2007), no. 6, 983--993.

\bibitem[L]{Ly} I.G. Lysenok, A set of defining relations for the Grigorchuk group, Mat. Zametki 38 (4) (1985) 503--516.

\bibitem[Mal]{Mal} A. I. Malcev, On isomorphic matrix representations of infinite groups. Rec. Math. [Mat. Sbornik] N.S.  8 (50),  (1940). 405--422.

\bibitem[Mar]{Mar} Gregory A. Margulis. Discrete subgroups of semisimple Lie groups, volume 17
of Ergebnisse der Mathematik und ihrer Grenzgebiete, Springer, 1991.

\bibitem[N]{N} P. M. Neumann, The $SQ$-universality of some finitely presented groups.
Collection of articles dedicated to the memory of Hanna Neumann, J. Austral. Math. Soc. 16 (1973), 1--6.

\bibitem[Ol91]{O1}
A.~Yu. Olshanskii, \newblock Periodic quotient groups of hyperbolic groups, \newblock { Mat. Sb.}, 182(4):543--567, 1991.

\bibitem[Ol93]{O2}
A.~Yu. Olshanskii, \newblock On residualing homomorphisms and {$G$}-subgroups of hyperbolic
  groups, \newblock { Internat. J. Algebra Comput.}, 3(4):365--409, 1993.



\bibitem[Ol95]{OlSq}
A.~Yu. Olshanskii,
\newblock {${\rm SQ}$}-universality of hyperbolic groups.
\newblock {\em Mat. Sb.}, 186(8):119--132, 1995.

\bibitem[Ol00]{OlBass} A. Yu. Olshanskii, On the Bass-Lubotzky question about quotients of hyperbolic groups.  J. Algebra  226  (2000),  no. 2, 807--817.


\bibitem[OO]{OO}
A.~Yu. Olshanskii and D.~V. Osin, \newblock Large groups and their periodic quotients, \newblock { Proc. Amer. Math. Soc.}, 136(3):753--759, 2008.



\bibitem[Pink]{Pink} Richard Pink. On the calculation of local terms in the
Lefschetz-Verdier trace formula and its application to a
conjecture of Deligne.  Ann. of Math. (2) 135 (1992), no. 3,
483--525.

\bibitem[Riv]{Rivin}
Igor Rivin.
\newblock Growth in free groups (and other stories).
\newblock {\em arXiv:math/9911076v2}, 1999.

\bibitem[SSch]{SSch} George S. Sacerdote, Paul E.  Schupp, SQ-universality in HNN groups and one relator groups.  J. London Math. Soc. (2)  7  (1974), 733--740.


\bibitem[San]{Sanov} I. N. Sanov. A property of a representation of a
free group. Doklady Akad. Nauk SSSR (N. S.) 57, (1947). 657--659.


\bibitem[SS]{SS} Mark Sapir, Iva \v Spakulov\' a, Almost all one-relator groups with at least three  generators are residually finite. preprint, arXiv math0809.4693, 2008.

\bibitem[SW]{SW} Mark Sapir, Daniel T. Wise. Ascending HNN extensions of
residually finite groups can be non-Hopfian and can have very few
finite quotients. J. Pure Appl. Algebra 166 (2002), no. 1-2,
191--202.

\bibitem[Wise1]{Wise0} Daniel T. Wise, A non-Hopfian automatic group.
J. Algebra 180 (1996), no. 3, 845--847.

\bibitem[Wise2]{Wise} Daniel T. Wise,  The residual finiteness of positive
one-relator groups. Comment. Math. Helv. 76 (2001), no. 2, 314--338.


\bibitem[Zel]{Z1} E.~I. Zelmanov, \newblock Solution of the restricted {B}urnside problem for groups of odd
  exponent, \newblock { Izv. Akad. Nauk SSSR Ser. Mat.}, 54(1):42--59, 221, 1990.


\end{thebibliography}
\end{document}